\documentclass{siamltex}
\usepackage{amsfonts,latexsym}
\usepackage{epsfig}
\usepackage{multirow}
\usepackage{amsmath}
\usepackage{graphicx}
\usepackage{caption}
\usepackage{subcaption}
\usepackage{amsmath}
\usepackage{color}
\usepackage{ulem}
\usepackage{algorithm,algorithmic}
\newcommand{\R}{\mathbb R}

\usepackage{epstopdf} 
\title{ A generalized matrix Krylov subspace method for TV regularization}
\author{A. H. Bentbib\thanks{Facult\'e des Sciences et Techniques-Gueliz, Laboratoire de Math\'ematiques Appliqu\'ees et Informatique, Morocco.. E-mail: {\tt a.bentbib@uca.ma}}
	\and 
	M. El Guide\thanks{Facult\'e des Sciences et Techniques-Gueliz, Laboratoire de Math\'ematiques Appliqu\'ees et Informatique, Morocco.. E-mail: {\tt mohamed.elguide@edu.uca.ac.ma}}
	\and 
	K. Jbilou\thanks{Universit\'e de Lille Nord de France, L.M.P.A, ULCO, 50 rue F. Buisson,
		BP699, F-62228 Calais-Cedex, France. E-mail: {\tt jbilou@univ-littoral.fr}}}
\date{}
 
\begin{document}
 \maketitle
 
\begin{abstract}
This paper presents an efficient algorithm to solve total variation (TV) regularizations of images contaminated by a both blur and  noise. The unconstrained structure of the problem suggests that one can solve a constrained optimization  problem by transforming the original unconstrained minimization problem to an equivalent constrained
minimization one. An augmented Lagrangian method is developed to handle the constraints when the model is given with matrix variables, and an alternating direction method
(ADM) is used to iteratively find solutions. The solutions of some sub-problems are belonging to subspaces generated by application of successive orthogonal projections onto a class of  generalized matrix Krylov subspaces of increasing dimension.
\end {abstract}


\section{Introduction}
In this paper we consider the solution of the following matrix equation
\begin{equation}\label{model1}
B=H_2XH_1^T,
\end{equation}
where $B$ is generally contaminated by noise. $H_1$ and $H_2$ are matrices of ill-determined rank, which makes the solution $X$ very
sensitive to perturbations in $B$. Discrete ill-posed problems of the form (\ref{model1}) arise, for instance, from the discretization of Fredholm
integral equations of the first kind in two space-dimensions,
\begin{equation}\label{fred2d}
\int\int_\Omega K(x,y,s,t)f(s,t)dsdt=g(x,y),\qquad (x,y)\in\Omega',
\end{equation}
where $\Omega$ and $\Omega'$  are rectangles in $\mathbb{R}^2$ and the kernel is separable
$$K(x,y,s,t)=k_1(x,s)k_2(y,t),\quad (x,y)\in\Omega',\quad (s,t)\in\Omega,$$
 The aim of this work is to solve this problem with application to one single channel and multichannel images. 
\subsection{Single channel images } 
For single channel images  we seek to recover an unknown vector from limited information. This problem is mathematically formulated as the following model
\begin{equation}\label{model2}
b=Hx,
\end{equation}
where $x\in\mathbb{R}^{mn}$ is a vector denoting the unknown solution, $b\in\mathbb{R}^{mn}$ is a vector denoting the observed data contaminated by noise and $H\in\mathbb{R}^{mn\times mn}$ is a linear map. The problem arises, for instance in image restoration \cite{An,BB,Ch,HNO,Ja}. In this paper we focus on the application to image restoration in which $x$ represents the unknown sharp image that is to be estimated from its blurry and noisy observation $b$. The matrix $H$ is the blurring operator characterized by a PSF describing this blur. Due to the ill-conditioning of the matrix $H$ and the presence of the noise, the problem (\ref{model2})  cannot be easily solved which means that the  minimization of only the fidelity
term  typically yields a meaningless computed solution.  Therefore, to stabilise  the recovered image, regularization is needed. There are several techniques to regularize the linear inverse problem given by
equation (\ref{model2}) ; see for example, \cite{GO,TA, ROF, VO}. All of these techniques stabilize the restoration process by adding a regularization term, depending on some a priori knowledge of  the unknown image, resulting in the model
\begin{equation}\label{model3}
\underset{x}{\text{min}}\{\|Hx-b\|^p_p+\mu\|\Phi(x)\|^q_q\},
\end{equation}
where $\Phi(x)$ is the regularizer that enforces the a priori knowledge and the parameter $\mu$ is used to balance the two terms.
This problem is referred to as $\ell_p-\ell_q$ minimization problem. Different choices of $\Phi(x)$, $p$ and $q$ lead to a wide variety of regularizers. Among them we find the well known Tikhonov regularization, where $\Phi$ is the identity
matrix, $p=2$ and $q=2$, see for example \cite{TA}. If the goal is to enforce sparsity on the solution, one can also consider $\Phi=I$, $p=2$ and $q=1$. Another well-known class of regularizers are based on total variation (TV), which is a better choice if the goal is to preserve sharp edges. In this case one let $\Phi$ to be the discrete gradient operator, see \cite{ROF}. The problem (\ref{model3}) has been studied in  many papers to  propose nonlinear optimization algorithms that can deal with the nonlinear properties of this  problem; see for example \cite{Vogel, WFN}. These techniques are computationally demanding if the main cost of computation  is the matrix-vector multiplication (MVM). It is our main goal to recover a good approximation of the unknown sharp image at  low computational cost. Because of some unique features in images, we seek an image
restoration algorithm that utilizes blur information, exploits
the spatially invariant properties. For this reason we suppose that the PSF is identical in all parts of the image and separates into horizontal and vertical components. Then the matrix $H$ is  the Kronecker product of two matrices $H_1$ and $H_2$,
\begin{equation}
H=H_1\otimes H_2=
\begin{bmatrix} h_{1,1}^{(1)} H_2 & h_{1,2}^{(1)} H_2 & \cdots & h_{1,n}^{(1)} H_2 \\ 
h_{2,1}^{(1)} H_2 & h_{2,2}^{(1)} H_2 & \cdots & h_{2,n}^{(1)} H_2\\
\vdots & \vdots & & \vdots \\
h_{n,1}^{(1)} H_2 & h_{n,2}^{(1)} H_2 & \cdots & h_{n,n}^{(1)} H_2 
\end{bmatrix}.
\end{equation}
In what follows we will need the \textsf{vec} and \textsf{mat} notations, which are a useful tools in transforming  the expression of matrix-vector product into a matrix-matrix product. Let the operator $\textsf{vec}$ transform a matrix $A=[a_{i,j}]\in\R^{m\times n}$ to a vector 
$a\in\R^{mn}$ by stacking the columns of $A$ from left to right, i.e,
\begin{equation}\label{avec}
a=[a_{1,1},a_{2,1},\ldots,a_{m,1},a_{1,2},a_{2,2},\ldots,a_{m,2},\ldots,a_{m,n}]^T,
\end{equation}
and let $\textsf{mat}$ be the inverse operator,  which transforms a vector (\ref{avec}) 
to an associated matrix $A=[a_{i,j}]\in\R^{m\times n}$. Thus,
\[
\textsf{vec}(A)=a,\qquad \textsf{mat}(a)=A.
\]
The Kronecker product satisfies the following relations for matrices $A,B,C,D,X$ of 
suitable sizes:
\begin{equation}\label{tensorrules}
\left.
\begin{array}{rcl}
(A\otimes B)\textsf{vec}(X)&=&\textsf{vec}(BXA^T),\\
(A\otimes B)^T&=&A^T\otimes B^T,\\
(AB)\otimes(CD)&=&(A\otimes C)(B\otimes D).
\end{array}
\right\}
\end{equation}
For $A,B\in\R^{m\times n}$, we define the inner product 
\begin{equation}\label{innerprod}
\langle A,B\rangle_F := \textsf{tr}(A^T B),
\end{equation}
where $\textsf{tr}(\cdot)$ denotes the trace. Notice  that
\begin{equation}\label{innerprodequiv}
\langle A,B\rangle_F = (\textsf{vec}(A))^T\textsf{vec}(B).
\end{equation}
The Frobenius norm is associated with this inner product,
\[
\|A\|_F := \langle A,A \rangle_F^{1/2},
\]
and it satisfies 
\begin{equation}\label{fnorm}
\|A\|_F=\|\textsf{vec}(A)\|_2.
\end{equation}
By using the properties (\ref{tensorrules}), the equation (\ref{model2}) can be rewritten as 
\begin{equation}
B=H_2XH_1^T,
\end{equation}
where $X=\textsf{mat}(x)$ and  $B=\textsf{mat}(b)$, which yields the model (\ref{model1}).
\subsection{Multichannel Images}
Recovering multichannel images from their blurry and noisy observations can be seen as a linear system of equations with multiple right-hand sides. The most commonly multichannel images is the RGB representation, which uses three channels; see \cite{GKCH, HNO}. It should be pointed out that the algorithms proposed in this paper
can be applied to the solution of Fredholm integral equations of the first kind in
two or more space dimensions and to the restoration of hyper-spectral images. The
latter kind of images generalize color images in that they allow more than three \textquotedblleft colors\textquotedblright; see, e.g., \cite{LNP}. If the channels are represented by $m\times n$ pixels, the full blurring model is described by the following form
\begin{equation}
b=Hx,
\end{equation}
where $b$ and $x$ in $ \mathbb{R}^{k mn}$, represent the blurred and noisy multichannel image and the original image  respectively. For an image with $k$ channels, they are given by
$$b=[b^{(1)}; b^{(2)};...;b^{(k)}],\quad x=[x^{(1)}; x^{(2)};...;x^{(k)}],$$
where $b^{(i)}$ and $x^{(i)}$  in $\mathbb{R}^{mn}$ are obtained  by stacking the columns of each channel on top of each other. The $kmn\times kmn$ multichannel blurring matrix $H$ is given by
\begin{equation}
H=H_1\otimes H_2,
\end{equation}
The matrix $H_2\in\mathbb{R}^{mn\times mn}$ represents  the same within-channel blurring in all the $k$ channels. The matrix $H_1$ of dimension $k\times k$ models the cross-channel blurring, which is the same for all pixels in the case of a spatially invariant blur. If $H_1=I$, the blurring is said to be within-channel. If no colour blurring arises (i.e., $H_1=I$), then $k$ independent deblurring problems are solved; hence the spatially invariant blurring model is given by
\begin{equation}\label{lsmrhs}
b_i=H_2x_i,\quad i=1,...,k.
\end{equation}
In this case, the goal is to model the blurring of $k$ channels image as a linear system of equations with $k$ right-hand sides. For this reason we let $B$ and $X$  in $\mathbb{R}^{mn\times k}$ to be denoted by $\left[b^{(1)}, b^{(2)},..., b^{(k)}\right]$ and $\left[x^{(1)}, x^{(2)},..., x^{(k)}\right]$, respectively. The optical blurring is then modeled by 
\begin{equation}\label{lsemrhs}
B=H_2X,
\end{equation}
which yields the model (\ref{model1}) with $H_1=I$.
When the spatially invariant cross-channel is present (i.e., $H_1\neq I$) and by using the Kronecker product properties, the following blurring model is to be  solved
\begin{equation}\label{linoper}
B=H_2XH_1^T,
\end{equation}
which also yields the model (\ref{model1}). Introduce the linear operator 
\begin{eqnarray*}
	\mathcal{H}:\mathbb{R}^{p\times q}& \rightarrow & \mathbb{R}^{p\times q} \\
	\mathcal{H}(X)&=&H_2XH_1^T.
\end{eqnarray*}
Its transpose is given by $\mathcal{H}^T(X) = H_2^TXH_1$. The problem (\ref{model1}) can be then expressed as 
\begin{equation*}
B=\mathcal{H}(X).
\end{equation*}

\noindent The total variation regularization  is known to be the most popular and effective techniques  for the images restoration. Given an image defined as a function $u:\Omega\longrightarrow \mathbb{R}$, where $\Omega$ is a bounded open subset of $\mathbb{R}^{2}$, the total variation (TV) of $u$ can be defined as
\begin{equation}\label{TV}
\text{TV}_{k}(u)=\int_{\Omega}\|\nabla u(x)\|_kdx,
\end{equation}
where $\nabla$ denotes  the gradient of $u$ and $\|.\|_k$ is a norm in $\mathbb{R}^{2}$. When $u$ is represented by $m\times n$ image $X$, a discrete form of (\ref{TV}) is always used, given by
\begin{equation}
\text{TV}_1(X)=\sum_{i=1}^{m}\sum_{j=1}^{n}\left(|\left(D_{1,n}X\right)_{ij}|+|\left(D_{1,m}X\right)_{ij}|\right)
\end{equation} 
in the anisotropic total variation case, or
\begin{equation}
\text{TV}_2(X)=\sum_{i=1}^{m}\sum_{j=1}^{n}\sqrt{\left(\left(D_{1,n}X\right)^2_{ij}+\left(D_{1,m}X\right)^2_{ij}\right)}
\end{equation}
in the isotropic total variation case.  $D_{1,m}$ and $D_{1,n}$ denote the finite
difference approximations of the horizontal and vertical first derivative operators, respectively, and they are defined as follows 
\begin{equation}
\begin{pmatrix}
D_{1,n}\\D_{1,m}
\end{pmatrix}X=\begin{pmatrix}
CX\\XC^T \end{pmatrix},
\end{equation}
where $$C:=\begin{bmatrix}
	-1&1&&\\
	&\ddots&\ddots&\\
	&&-1&1
\end{bmatrix}\in\mathbb{R}^{d-1\times d},$$
where  $d$ is the number of pixels in each row
 and column of the image considered. For the ill-posed image restoration problem \eqref{model1}, the resulting matrices $H_1$ and $H_2$ are ill-conditioned.
By regularization of the problem (\ref{model1}), we solve as a special case  one of the following  matrix problems:
\begin{equation}
\underset{X}{ \text{min}}\left(\|\mathcal{H}(X)-B\|_F^2+\mu\text{TV}_k(X)\right),\quad k=1,2\label{TVL2}
\end{equation}
or
\begin{equation}
\underset{X}{ \text{min}}\left(\|\mathcal{H}(X)-B\|_{1,1}+\mu\text{TV}_k(X)\right)\quad k=1,2.\label{TVL1}
\end{equation}
where $\|.\|_{1,1}$ is the $\ell_1$ norm and $\mu$ is a regularization parameter.
 Problems (\ref{TVL2}) and (\ref{TVL1}) are refereed to as TV/L2 and TV/L1 minimization, respectively.
\section{TV/L2 minimization problem}
In this section we consider the solution of the following TV/L2 minimization problem
\begin{equation}\label{TVL2prob}
\underset{X}{ \text{min}}\left(\|\mathcal{H}(X)-B\|_F^2+\mu\text{TV}_2(X)\right).
\end{equation}
The model (\ref{TVL2prob}) is very difficult to solve directly due to the non-differentiability and non-linearity of the TV term. It is our goal to develop an efficient TV minimization scheme to handle this problem. The core idea is based on augmented Lagrangian method (ALM) \cite{Hest,Pow} and alternating direction method (ADM) \cite{GM}. The idea of ALM is to transform  the unconstrained minimization task (\ref{TVL2prob}) into an  equivalent constrained optimization problem, and then add a quadratic penalty term instead of the constraint violation with the multipliers. The idea of ADM is to decompose the transformed  minimization problem
into three easier and smaller sub-problems such
that some involved variables can be minimized
separately and alternatively. Let us begin by  considering the equivalent equality-constrained problem of (\ref{TVL2prob}). We first notice that the minimization problem (\ref{TVL2prob}) can be rewritten as 
\begin{equation}\label{TVL2cons}
\underset{X,M^{(n)},M^{(m)}}{ \text{min}}\left(\|\mathcal{H}(X)-B\|_F^2+\mu\sum_{i=1}^{m}\sum_{j=1}^{n}\|M_{i,j}\|_2\right),
\end{equation}
$$\text{subject to}\quad D_{1,n}X=M^{(n)},\quad D_{1,m}X=M^{(m)}.$$
where $M_{i,j}=\left[\left(D_{1,n}X\right)_{ij},\left(D_{1,m}X\right)_{ij}\right]$. If we set $M_{i,j}^{(n)}=\left(D_{1,n}X\right)_{ij}$ and $M_{i,j}^{(m)}=\left(D_{1,m}X\right)_{ij}$
This constrained problem can be also formulated as
\begin{eqnarray}\label{cop}
&&\underset{}{ \text{min}}\quad F(X) +G(Y),\\\nonumber
&&\text{subject to}\quad DX=Y,
\end{eqnarray} 
where,
$$F(X)=\|\mathcal{H}(X)-B\|_F^2,\quad G(Y)=\mu\sum_{i=1}^{m}\sum_{j=1}^{n}\|M_{i,j}\|_2, \quad D=\begin{pmatrix}
D_{1,n}\\D_{1,m}
\end{pmatrix},\quad Y=\begin{pmatrix}
M^{(n)}\\M^{(m)} \end{pmatrix}$$

The augmented Lagrangian
function of (\ref{cop}) is defined as
\begin{equation}\label{TVL2lagrangian}
\mathcal{L}_{\beta}\left(X,Y,Z\right)=F(X)+G(Y)+\left<DX-Y,Z\right>+\frac{\beta}{2}\|DX-Y\|_F^2,
\end{equation}
where $Z\in\mathbb{R}^{2m\times n}$ 
is the Lagrange multiplier of the linear
constraint and $\beta>0$ is the penalty parameter for the violation of this linear constraint.\\
To solve the nonlinear problem (\ref{TVL2prob}), we find the saddle point of the Lagrangian (\ref{TVL2lagrangian}) by using the ADM method. The idea of this method is to apply an alternating minimization iterative procedure, namely, for $k = 0, 1, . . . ,$ we solve 
\begin{eqnarray}
(X_{k+1},Y_{k+1})=\underset{X}{\text{arg min }} \mathcal{L}_{\beta}(X,Y,Z_k)\label{TVL2Xprob}.
\end{eqnarray}
The Lagrange multiplier is updated by
\begin{eqnarray}\label{ZTVL2}
Z_{k+1}&=&Z_k+\beta\left(DX_{k+1}-Y_{k+1}\right).
\end{eqnarray} 
\subsection{Solving the Y-problem} Given $X$, 
$Y_{k+1}$ can be obtained by solving 
\begin{equation}\label{YTVL2}
\underset{Y}{\text{ min }} \mu\sum_{i=1}^{m}\sum_{j=1}^{n}\|M_{i,j}\|_2+\frac{\beta}{2}\|DX-Y\|_F^2+\left<DX-Y,Z_k\right>_F
\end{equation}
which is equivalent to solve
\begin{equation}
\underset{Y}{\text{ min }} \mu\sum_{i=1}^{m}\sum_{j=1}^{n}\|M_{i,j}\|_2+\frac{\beta}{2}\left\|\begin{pmatrix}
M^{(n)}\\M^{(m)} \end{pmatrix}-\begin{pmatrix}
D_{1,n}X\\D_{1,m}X
\end{pmatrix}-\frac{1}{\beta}\begin{pmatrix}
Z^{(1)}_k\\Z^{(2)}_k
\end{pmatrix}\right\|_F^2
\end{equation}
which is also equivalent to solve the so-called M-subproblem
\begin{equation}\label{TVL1M}
\underset{M_{i,j}}{\text{ min }} \sum_{i=1}^{m}\sum_{j=1}^{n}\mu\|M_{i,j}\|_2+\frac{\beta}{2}\left|
M^{(n)}_{ij}-
K_{ij}\right|_F^2+\frac{\beta}{2}\left|
M^{(m)}_{ij}-
L_{ij}\right|_F^2
\end{equation}
where $K_{ij}=\left(D_{1,n}X\right)_{ij}
+\frac{1}{\beta}\left(Z^{(1)}_k\right)_{ij}$ and $L_{ij}=\left(D_{1,m}X\right)_{ij}
+\frac{1}{\beta}
\left(Z^{(2)}_k\right)_{ij}$. 
To solve (\ref{TVL1M}) we use following  well-known two dimensional shrinkage formula \cite{Li}
\begin{equation}
\textbf{Shrink}(y,\gamma,\delta)=\text{max}\left\{\left\|y+\frac{\gamma}{\delta}\right\|_2-\frac{1}{\delta},0\right\}\frac{y+\gamma/\delta}{\|y+\gamma/\delta\|_2},
\end{equation}
where the convention 0·(0/0) = 0 is followed. The solution of (\ref{TVL1M}) is then given by
\begin{equation}\label{Ysol}
M_{i,j}=\text{max}\left\{\|T_{i,j}\|_2-\frac{\mu}{\beta},0\right\}\frac{T_{i,j}}{\|T_{i,j}\|_2},
\end{equation}
where $T_{i,j}=\left[\left(D_{1,n}X_{k}\right)_{i,j}
+\frac{1}{\beta}
\left(Z^{(1)}_k\right)_{i,j},\left(D_{1,m}X_{k}\right)_{i,j}
+\frac{1}{\beta}
\left(Z^{(2)}_k\right)_{i,j}\right].$\\
For the anisotropic case we solve the following problem
\begin{equation}
\underset{M_{i,j}}{\text{ min }} \sum_{i=1}^{m}\sum_{j=1}^{n}\mu\|M_{i,j}\|_1+\frac{\beta}{2}\left|
M^{(n)}_{ij}-
K_{ij}\right|_F^2+\frac{\beta}{2}\left|
M^{(m)}_{ij}-
L_{ij}\right|_F^2
\end{equation}
which can be also solved by the one dimensional shrinkage formula. This gives
\begin{eqnarray}
M^{(n)}_{ij}&=&\text{max}\left\{K_{ij}-\frac{\mu}{\beta},0\right\}.\text{sign}\left(K_{ij}\right),\\
M^{(m)}_{ij}&=&\text{max}\left\{L_{ij}-\frac{\mu}{\beta},0\right\}.\text{sign}\left(L_{ij}\right),
\end{eqnarray}
\subsection{Solving the X-problem} Given $Y$, 
$X_{k+1}$ can be obtained by solving 
\begin{equation}\label{XTVL2}
\underset{X}{\text{min }}\frac{\beta}{2}\|DX-Y\|_F^2+\left<DX-Y,Z_k\right>_F+\|\mathcal{H}(X)-B\|_F^2.
\end{equation}
This problem can be also solved by considering the following normal equation
\begin{equation}
H_1^TH_1XH_2^TH_2+\beta D^TDX=H_1^TBH_2+ D^T\left(\beta Y-Z_k\right).
\end{equation}
The linear matrix equation can be rewritten in the following form
\begin{equation}\label{TVL2sylv}
A_1XA_2+A_3XA_4=E_k,\quad k=1,...,
\end{equation}
where $A_1=H_1^TH_1$, $A_2=H_2^TH_2$, $A_3=\beta D^TD$, $A_4=I$ and $E_k=H_1^TBH_2+D^T\left(\beta Y-Z_k\right)$. The equation (\ref{TVL2sylv}) is refereed to as the generalized Sylvester matrix equation.   We will see in section \ref{sec4}  how to compute approximate solutions to those matrix equations 
\subsection{Convergence analysis of TV/L2 problem}
For the vector case, many convergence results have been proposed in the literature ; see for instance \cite{Glow, HLHY}. For completeness,
we give a proof here for the matrix case. A function $\Psi$ is said to be proper if the domain of $\Psi$ denoted by $\textbf{dom}\Psi:=\left\{U\in \mathbb{R}^{p\times q}, \Psi(U)<\infty\right\}$ is not empty. For the problem (\ref{cop}), $F$ and $G$ are closed proper convex functions.  According to \cite{FP, Rock}, the problem (\ref{cop}) is solvable, i.e., there exist $X_*$ and $Y_*$, not necessarily unique that minimize (\ref{cop}).  Let $\mathcal{W}=\Omega\times \mathcal{Y}\times \mathbb{R}^{p\times q}$, where $\Omega$ and $\mathcal{Y}$ are given closed and convex nonempty  sets. The saddle-point problem is equivalent to finding $(X_*,Y_*,Z_*)\in\mathcal{W}$ such that
\begin{equation}
\mathcal{L}_{\beta}(X_*,Y_*,Z)\leq\mathcal{L}_{\beta}(X_*,Y_*,Z_*)\leq\mathcal{L}_{\beta}(X,Y,Z_*),\qquad \forall \left(X,Y,Z\right)\in\mathcal{W}.
\end{equation}
The properties of the relation between the saddle-points of $\mathcal{L}_{\beta}$ and  $\mathcal{L}_{0}$ and the solution of (\ref{cop}) are stated by the following theorem from \cite{Glow}

\begin{theorem}\label{Thm1}
	$(X_*,Y_*,Z_*)$ is a saddle-point of $\mathcal{L}_{0}$ if and only if $(X_*,Y_*,Z_*)$ is a saddle-point of $\mathcal{L}_{\beta}$ $\forall \beta >0$. Moreover $(X_*,Y_*)$ is a solution of (\ref{cop}).
\end{theorem}\\

We will see in what follows how this theorem can be used to give the convergence of $\left(X_{k+1},Y_{k+1}\right)$. It should be pointed out that the idea of our proof follows the convergence results in \cite{BPCPE}.
\begin{theorem}
	Assume that $(X_*,Y_*,Z_*)$ is a saddle-point of $\mathcal{L}_{\beta}$ $\forall \beta >0$. The sequence $(X_{k+1},Y_{k+1},Z_{k+1})$ generated by Algorithm 1 satisfies

	\begin{enumerate}
		\item $\lim\limits_{k \rightarrow +\infty} F(X_{k+1}) + G(Y_{k+1})=F(X_*)+G(Y_*)$,
		\item  $\lim\limits_{k \rightarrow +\infty}\|DX_{k+1}-Y_{k+1}\|_F=0,$
	\end{enumerate}
\end{theorem}

$\textbf{Proof}$ \quad In order to show the convergence of this theorem, it suffice to show that the  non-negative function 
\begin{equation}
F^k=\frac{1}{\beta}\|Z_k-Z_*\|_F^2+\beta\|X_k-X_*\|_F^2
\end{equation}
decreases at each iteration. Let us define $S_k$, $M_k$ and $M_*$ as
$$S_k=DX_{k}-Y_{k}, \quad M_k=F(X_{k}) + G(Y_{k}), \quad M_*=F(X_*)+G(Y_*).$$
In the following we show
\begin{equation}\label{ineq}
F^{k+1}\leq F^k-\beta \|S_{k+1}\|_F^2-\beta\|Y_{k+1}-Y_k\|_F^2.
\end{equation}
Since $(X_*,Y_*,Z_*)$ is a saddle-point of $\mathcal{L}_{\beta}$ $\forall \beta >0$, it follows from Theorem \ref{Thm1} that $(X_*,Y_*,Z_*)$ is also  a saddle-point of $\mathcal{L}_{0}$. This is characterized by
\begin{equation}\label{unL}
\mathcal{L}_{0}(X_*,Y_*,Z)\leq\mathcal{L}_{0}(X_*,Y_*,Z_*)\leq\mathcal{L}_{0}(X,Y,Z_*),\qquad \forall \left(X,Y,Z\right)\in\mathcal{W}.
\end{equation}
From the second inequality of (\ref{unL}), we have
\begin{equation}\label{q1}
M_*-M_{k+1}\leq \left< S_{k+1},Z_*\right>_F.
\end{equation}
In the oder hand, $X_{k+1}$ is a minimizer of $\mathcal{L}_{\beta}$ $\forall \beta >0$, this implies that the optimality conditions reads 
\begin{equation}
2\mathcal{H}^T\left(\mathcal{H}(X_{k+1})-B\right)+D^T(Z_k+\beta(DX_{k+1}-Y_k))=0.
\end{equation}
By plugging $Z_k=Z_{k+1}-\beta(DX_{k+1}-Y_{k+1})$ and rearranging we obtain
\begin{equation}
2\mathcal{H}^T\left(\mathcal{H}(X_{k+1})-B\right)+D^T(Z_{k+1}-\beta(Y_{k+1}-Y_k))=0,
\end{equation}
which means that $X_{k+1}$ minimizes
\begin{equation}
F(X)+\left<Z_{k+1}+\beta(Y_{k+1}-Y_k),DX\right>_F.
\end{equation}
It follows that 
\begin{equation}\label{c1}
F(X_{k+1})-F(X_*)\leq\left<Z_{k+1}+\beta(Y_{k+1}-Y_k),DX_*\right>_F-\left<Z_{k+1}+\beta(Y_{k+1}-Y_k),DX_{k+1}\right>_F.
\end{equation}
A similar argument shows that
\begin{equation}\label{c2}
G(Y_{k+1})-G(Y_*)\leq\left<Z_{k+1},Y_{k+1}\right>-\left<Z_{k+1},Y_*\right>_F.
\end{equation}
Adding (\ref{c1}) and (\ref{c2}) and using $DX_*=Y_*$ implies
\begin{equation}\label{q2}
M_{k+1}-M_*\leq -\left<S_{k+1},Y_{k+1}\right>_F-\left<\beta(Y_{k+1}-Y_k),S_{k+1}+(Y_{k+1}-Y_*)\right>_F.
\end{equation}
Adding (\ref{q1}) and (\ref{q2}) and multiplying through by
2  gives
\begin{equation}\label{Y1}
2\left<S_{k+1},Z_{k+1}-Z_*\right>_F+2\left<\beta(Y_{k+1}-Y_k),S_{k+1}\right>_F+2\left<\beta(Y_{k+1}-Y_k),(Y_{k+1}-Y_*)\right>_F\leq 0
\end{equation}
The inequality (\ref{ineq}) will hold by rewriting each term of the inequality (\ref{Y1}). Let us begin with its first term. 
Substituting $Z_{k+1}=Z_k+\beta S_{k+1}$ gives
\begin{equation}\label{d1}
2\left<S_{k+1},Z_{k+1}-Z_k\right>_F=2\left<S_{k+1},Z_{k}-Z_*\right>_F+\beta\|S_{k+1}\|_F^2+\beta\|S_{k+1}\|_F^2.
\end{equation}
Since $S_{k+1}=\frac{1}{\beta}\left(Z_{k+1}-Z_k\right)$, it follows that the first two terms of the right hand side of (\ref{d1}) can be written as
\begin{equation}\label{k1}
\frac{2}{\beta}\left<Z_{k+1}-Z_k,Z_{k}-Z_*\right>_F+\frac{1}{\beta}\|Z_{k+1}-Z_k\|_F^2.
\end{equation}
Substituting $Z_{k+1}-Z_k=(Z_{k+1}-Z_*)-(Z_k-Z_*)$, shows that (\ref{k1}) can be written as 
\begin{equation}\label{u1}
\frac{1}{\beta}\left(\|Z_{k+1}-Z_*\|_F^2-\|Z_k-Z_*\|_F^2\right).
\end{equation}
We turn now to the remaining terms, i.e.,
\begin{equation}\label{remter}
\beta\|S_{k+1}\|_F^2+2\left<\beta(Y_{k+1}-Y_k),S_{k+1}\right>_F+2\left<\beta(Y_{k+1}-Y_k),(Y_{k+1}-Y_*)\right>_F
\end{equation}
Substituting $Y_{k+1}-Y_*=(Y_{k+1}-Y_k)+(Y_{k+1}-Y_*)$ shows that (\ref{remter}) can be expressed as 
\begin{equation}\label{m1}
\beta\|S_{k+1}+(Y_{k+1}-Y_k)\|_F^2+\beta\|Y_{k+1}-Y_k\|_F^2+2\beta\left<Y_{k+1}-Y_k,Y_k-Y_*\right>_F.
\end{equation}
Substituting $Y_{k+1}-Y_k=(Y_{k+1}-Y_*)-(Y_k-Y_*)$ in  the last two terms shows that (\ref{m1}) can be expressed as
\begin{equation}\label{u2}
\beta\|S_{k+1}+(Y_{k+1}-Y_k)\|_F^2+\beta\left(\|Y_{k+1}-Y_*\|_F^2-\|Y_{k}-Y_*\|_F^2\right)
\end{equation}
Using (\ref{u1}) and (\ref{u2}) shows that (\ref{Y1}) can be expressed as
\begin{equation}
F^k-F^{k+1}\geq \beta \|S_{k+1}+(Y_{k+1}-Y_k)\|_F^2.
\end{equation}
To show (\ref{ineq}), it is now suffice to show that $2\beta\left<S_{k+1},Y_{k+1}-Y_k\right>_F\geq0$. Since $(X_k,Y_k,Z_k)$ and $(X_{k+1},Y_{k+1},Z_{k+1})$  are also  minimizers  of $\mathcal{L}_\beta$, we have as in (\ref{c2})
\begin{equation}\label{s1}
G(Y_{k+1})-G(Y_k)\leq\left<Z_{k+1},Y_{k+1}\right>_F-\left<Z_{k+1},Y_k\right>_F,
\end{equation}
and 
\begin{equation}\label{s2}
G(Y_{k})-G(Y_{k+1})\leq\left<Z_{k},Y_{k}\right>_F-\left<Z_{k},Y_{k+1}\right>_F.
\end{equation}
It follows by addition of (\ref{s1}) and (\ref{s2}) that,
\begin{equation}
\left<Y_{k+1}-Y_k,Z_{k+1}-Z_k\right>\geq0.
\end{equation}
Substituting  $Z_{k+1}-Z_k=\beta S_{k+1}$ shows that $2\beta\left<S_{k+1},Y_{k+1}-Y_k\right>\geq0$.
From (\ref{ineq}) it follows that 
\begin{equation}
\beta\sum_{k=0}^{\infty}\left(\|S_{k+1}\|_F^2-\beta\|Y_{k+1}-Y_k\|_F^2\right)\leq F^0,
\end{equation}
which implies that $S_{k+1}\longrightarrow 0$ and $Y_{k+1}-Y_k\longrightarrow0$ as $k\longrightarrow \infty$. It follows then from (\ref{q1}) and (\ref{q2}) that $\lim\limits_{k \rightarrow +\infty} F(X_{k+1}) + G(X_{k+1})=F(X_*)+G(X_*)$,
\section{TV/L1 minimization problem}
In this section we consider the following regularized minimization problem
\begin{equation}\label{TVL1model}
\underset{X}{ \text{min}}{\|\mathcal{H}(X)-B\|_{1,1}+\mu\text{TV}_2(X)}
\end{equation}
 We first notice that the minimization problem (\ref{TVL1model}) can be rewritten as
\begin{equation}
\underset{X}{ \text{min}}\left(\|\mathcal{H}(X)-B\|_{1,1}+\mu\sum_{i=1}^{m}\sum_{j=1}^{n}\|M_{i,j}\|_2\right),
\end{equation}
 then,
the constraint violation of the problem (\ref{TVL1model})  can be written as follows
\begin{equation}
\underset{X,R,M^{(n)},M^{(m)}}{ \text{min}}\left(\|R-B\|_{1,1}+\mu\sum_{i=1}^{m}\sum_{j=1}^{n}\|M_{i,j}\|_2\right),
\end{equation}
$$\text{subject to}\quad D_{1,n}X=M^{(n)},\quad D_{1,m}X=M^{(m)},\quad R=\mathcal{H}(X).$$
This constrained  problem can be also reformulated as
\begin{eqnarray}\label{constraint}
&&\underset{}{ \text{min}}\quad F(R) +G(Y),\\\nonumber
&&\text{subject to}\quad DX=Y, \quad \mathcal{H}(X)=R
\end{eqnarray} 
where,
$$F(R)=\|R-B\|_{1,1},\quad G(Y)=\mu\sum_{i=1}^{m}\sum_{j=1}^{n}\|M_{i,j}\|_2, \quad D=\begin{pmatrix}
D_{1,n}\\D_{1,m}
\end{pmatrix},\quad Y=\begin{pmatrix}
M^{(n)}\\M^{(m)} \end{pmatrix}, $$
The problem now fits the framework of the augmented Lagrangian method \cite{Hest,Pow} which puts a quadratic penalty term instead of the constraint in
the objective function and introducing explicit Lagrangian multipliers at each iteration into the objective function.  The augmented Lagrangian
function of (\ref{constraint}) is defined as follows
\begin{eqnarray}\label{TVL1lagrangian}
&&\mathcal{L}\left(X,R,Y,Z,W\right) =\\\nonumber&&F(R) +  G(Y)+\frac{\beta}{2}\|DX-Y\|_F^2+\left<DX-Y,Z\right>_F+\frac{\rho}{2}\|\mathcal{H}(X)-R\|_F^2+\left<\mathcal{H}(X)-R,W\right>_F
\end{eqnarray}

$Z\in\mathbb{R}^{2m\times n}$ and $W\in\mathbb{R}^{m\times n}$ are the Lagrange multipliers of the linear
constraint $DX=Y$ and $R=\mathcal{H}(X)$ , respectively. The parameters $\beta>0$ and $\rho>0$ are the penalty parameters for the violation of the linear constraint. \\
Again, we use the ADM method to solve the nonlinear problem (\ref{TVL1model}), by finding  the saddle point of the Lagrangian (\ref{TVL1lagrangian}). Therefore, for $k=0,1,...$ we solve 
\begin{equation}
\left(X_{k},R_{k},Y_k\right)=\underset{X,R,Y}{\text{arg min }} \mathcal{L}_{\beta,\rho}(X,R,Y,Z_k,W_k).\label{TVL1Xprob}
\end{equation}
The Lagrange multipliers are updated by
\begin{eqnarray}\label{lgup}
Z_{k+1}&=&Z_k+\beta\left(DX_{k}-Y_{k}\right).\nonumber\\
W_{k+1}&=&W_k+\rho\left(\mathcal{H}(X_{k})-R_{k}\right).
\end{eqnarray} 
Next, we will see how to solve the problems \eqref{TVL1Xprob}, to determine the iterates $X_{k}$, $Y_{k}$ and $R_{k}$
\subsection{Solving the X-problem} Given $Y$ and $R$, $X_k$ can be obtained by solving the minimization problem
\begin{eqnarray}\label{XTVL1diff}
\underset{X}{\text{min }} \frac{\beta}{2}\|DX-Y\|_F^2+\left<DX-Y,Z_k\right>_F+\frac{\rho}{2}\|\mathcal{H}(X)-R\|_F^2+\left<\mathcal{H}(X)-R,W_k\right>_F
\end{eqnarray}
The problem \eqref{XTVL1diff} is  now  continuously differentiable at $X$. Therefore, it can be solved by considering the following normal equation
\begin{equation}\label{lin1}
\rho H_1^TH_1XH_2^TH_2+\beta D^TDX= H_1^T\left(\rho R-W_k\right)H_2+D^T\left(\beta Y-Z_k\right).
\end{equation}
The linear matrix equation \eqref{lin1} can be rewritten in the following form
\begin{equation}\label{TVL1Sylv}
A_1XA_2+A_3XA_4=E_k,
\end{equation}
where $A_1=\rho H_1^TH_1$, $A_2=H_2^TH_2$, $A_3=\beta D^TD$, $A_4=I$ and $E_k= H_1^T\left(\rho R-W_k\right)H_2$\\$+D^T\left(\beta Y-Z_k\right)$.\\ 
 The equation (\ref{TVL1Sylv}) is refereed to as the generalized Sylvester matrix equation.
\subsection{Solving the R-problem} Given $X$, the iterate $R_{k}$ can be obtained by solving the minimization problem 
\begin{equation}\label{Rprob}
\underset{R}{\text{ min }} \|R-B\|_{1,1}+\frac{\rho}{2}\left\|\mathcal{H}(X)-R\right\|_F^2+\left<\mathcal{H}(X)-R,W_k\right>_F.
\end{equation}
Therefore, by using the following well-known one-dimensional Shrinkage formula \cite{Li}
\begin{equation}
\textbf{Shrink}(y,\gamma,\delta)=\text{max}\left\{\left|y+\frac{\gamma}{\delta}\right|-\frac{1}{\delta},0\right\}.\text{sign}\left(y+\frac{\gamma}{\delta}\right),
\end{equation}
the  minimizer of (\ref{Rprob}) is then given by
\begin{equation}\label{RTVL1}
\text{max}\left\{\left|\mathcal{H}(X)-B+\frac{1}{\rho}W\right|-\frac{1}{\rho},0\right\}.\text{sign}\left(\mathcal{H}(X)-B+\frac{1}{\rho}W\right)
\end{equation}
\subsection{Solving the Y-problem}
Given $X$ and $R$, we compute the iterates $Y_{k}$  by solving  the problem
\begin{equation}\label{YTVL1}
\underset{Y}{\text{ min }} \mu\sum_{i=1}^{m}\sum_{j=1}^{n}\|M_{i,j}\|_2+\frac{\beta}{2}\|DX-Y\|_F^2+\left<DX-Y,Z_k\right>_F
\end{equation}
This solution can be obtained by equation (\ref{Ysol}), since the minimization problem (\ref{YTVL1})  is the same
as that of TV/L2. 
\subsection{Convergence analysis of TV/L1 problem}
In this subsection we study the convergence of Algorithm 2 used to solve the TV/L1 problem. Note that the convergence study for TV/L2 does not hold for TV/L1 problem since in general $\beta\neq \rho$ in (\ref{TVL1lagrangian}). For the problem (\ref{constraint}), $F$ and $G$ are closed proper convex functions.  According to \cite{FP, Rock}, the problem (\ref{constraint}) is solvable, i.e., there exist $R^*$ and $Y^*$, not necessarily unique that minimize (\ref{constraint}).  Let $\mathcal{W}=\Omega\times \mathcal{Y}\times \mathcal{X}\times\mathbb{R}^{2m\times n}\times\mathbb{R}^{m\times n}$, where $\Omega$, $\mathcal{X}$ and $\mathcal{Y}$ are given closed and convex nonempty  sets. The saddle-point problem is equivalent to finding $(X_*,R_*,Y_*,Z_*,W_*)\in\mathcal{W}$ such that
\begin{eqnarray}\label{sadlpt}
&\mathcal{L}_{\beta,\rho}&(X_*,R_*,Y_*,Z,W)\leq\mathcal{L}_{\beta,\rho}(X_*,R_*,Y_*,Z_*,W_*)\leq\mathcal{L}_{\beta,\rho}(X,R,Y,Z_*,W_*),\nonumber \\&&\forall \left(X,R,Y,Z,W\right)\in\mathcal{W}.
\end{eqnarray}
The properties of the relation between the saddle-points of $\mathcal{L}_{\beta,\rho}$ and  the solution of (\ref{constraint}) are stated by the following theorem from \cite{WZT}
\begin{theorem}\label{Thm3}
$X_*$ is a solution of (\ref{TVL1model}) if and only if there exist $(R_*,Y_*)\in\mathcal{Y}\times \mathcal{X}$ and $(Y_*,Z_*)\in \mathbb{R}^{2m\times n}\times\mathbb{R}^{m\times n}$ such that $(X_*,R_*,Y_*,Z_*,W_*)$ is a saddle-point of (\ref{sadlpt})
\end{theorem}\\
The convergence of ADM for TV/L1 has been well studied in the literature in the context of vectors; see, e.g., \cite{WZT}. Our TV/L1 problem is  a model with matrix variables, it is our aim to give a similar convergence results for the matrix case
\begin{theorem}
	Assume that $(X_*,R_*,Y_*,Z_*,W_*)$ is a saddle-point of $\mathcal{L}_{\beta,\rho}$. The sequence $(X_{k},R_{k},Y_{k},Z_{k},W_{k})$ generated by Algorithm 2 satisfies
	
	\begin{enumerate}
		\item $\lim\limits_{k \rightarrow +\infty} F(R_{k}) + G(Y_{k})=F(R_*)+G(Y_*)$,
		\item  $\lim\limits_{k \rightarrow +\infty}\|DX_{k}-Y_{k}\|_F=0,$
		\item  $\lim\limits_{k \rightarrow +\infty}\|\mathcal{H}(X_{k})-R_{k}\|_F=0.$
	\end{enumerate}
\end{theorem}
$\textbf{Proof}$ \quad 
From the first inequality of (\ref{sadlpt}) it follows that $\forall (Z,W)\in \mathbb{R}^{2m\times n}\times\mathbb{R}^{m\times n}$
\begin{equation}
\left<DX_*-Y_*,Z_*\right>_F+\left<\mathcal{H}(X_*)-R_*,W_*\right>_F\leq\left<DX_*-Y_*,Z\right>_F+\left<\mathcal{H}(X_*)-R_*,W\right>_F,
\end{equation}
which obviously implies that
\begin{eqnarray}\label{rls}
DX_*=Y_*,\nonumber\\
\mathcal{H}(X_*)=R_*.
\end{eqnarray}
Let us define the following quantities
$$\overline{Z}_k=Z_k-Z_*,\quad\overline{W}_k=W_k-W_*,\quad\overline{X}_{k}=X_{k}-X_*,\quad\overline{R}_{k}=R_{k}-R_*,\quad\overline{Y}_{k}=Y_{k}-Y_*. $$
With the relationship (\ref{rls}) together with (\ref{lgup}), we can define
\begin{eqnarray}\label{lagup}
\overline{Z}_{k+1}&=&\overline{Z}_k+\beta\left(D\overline{X}_{k}-\overline{Y}_{k}\right)\\
\overline{W}_{k+1}&=&\overline{W}_k+\rho\left(\mathcal{H}(\overline{X}_{k})-\overline{R}_{k}\right)
\end{eqnarray}
In order to show the convergence, it suffice to show that $\left(\beta\|\overline{Z}_k\|_F^2+\rho\|\overline{W}_k\|_F^2\right)$ decreases at each iteration. In the following we show that
\begin{eqnarray}
\left(\beta\|\overline{Z}_k\|_F^2+\rho\|\overline{W}_k\|_F^2\right)&-&\left(\beta\|\overline{Z}_{k+1}\|_F^2+\rho\|\overline{W}_{k+1}\|_F^2\right)\\&\geq&\beta^2\rho\|D\overline{X}_{k}-\overline{Y}_{k}\|_F^2+\beta\rho^2\|\mathcal{H}(\overline{X}_{k})-\overline{R}_{k}\|_F^2.
\end{eqnarray}
For $\left(X,R,Y\right)=\left(X_{k},R_{k},Y_{k}\right)$ in (\ref{sadlpt}) , the second equality implies
\begin{eqnarray}
&&\left<D^TZ_*,X_{k}-X_*\right>_F+\beta\left<D^T(Y_*-DX_*),X_{k}-X_*\right>_F\\\nonumber&+&\left<W_*,-\mathcal{H}(X_{k}-X_*)\right>_F+\rho\left<Z_*-\mathcal{H}(X_*),-\mathcal{H}(X_{k}-X_*)\right>_F\geq0,
\end{eqnarray}
\begin{equation}
F(R_{k})-F(R_*)+\left<W_*,R_{k}-R_*\right>_F+\rho\left<R_*-\mathcal{H}(X_*),R_{k}-R_*)\right>_F\geq0,
\end{equation}
\begin{equation}
G(Y_{k})-G(Y_*)+\left<Z_*,Y_{k}-Y_*\right>_F+\beta\left<Y_*-DX_*,Y_{k}-Y_*)\right>_F\geq0.
\end{equation}
Since $\left(X_{k},R_{k},Y_{k}\right)$ is also a saddle-point of $\mathcal{L}_{\beta,\rho}$, for $\left(X,R,Y\right)=\left(X_*,R_*,Y_*\right)$ the second equality of (\ref{sadlpt})  implies
\begin{eqnarray}\label{1ineq}
&&\left<D^TZ_{k},X_*-X_{k}\right>_F+\beta\left<D^T(Y_{k}-DX_{k}),X_*-X_{k}\right>_F\\\nonumber&+&\left<W_{k},-\mathcal{H}(X_*-X_{k})\right>_F+\rho\left<Z_{k}-\mathcal{H}(X_{k}),-\mathcal{H}(X_*-X_{k})\right>_F\geq0,
\end{eqnarray}
\begin{equation}\label{2ineq}
F(R_*)-F(R_{k})+\left<W_{k},R_*-R_{k}\right>_F+\rho\left<R_{k}-\mathcal{H}(X_{k}),R_*-R_{k})\right>\geq0,
\end{equation}
\begin{equation}\label{3ineq}
G(Y_*)-G(Y_{k})+\left<Z_{k},Y_*-Y_{k}\right>_F+\beta\left<Y_{k}-DX_{k},Y_*-Y_{k})\right>\geq0.
\end{equation}
By addition , regrouping terms, and multiplying through by
$\beta\rho$ gives
\begin{equation}\label{keyf}
-\beta\rho \left<\overline{Z}_k,D\overline{X}_{k}-\overline{Y}_{k}\right>-\beta\rho \left<\overline{W}_k,\mathcal{H}(\overline{X}_{k})-\overline{Z}_{k}\right>\geq\beta^2\rho\|D\overline{X}_{k}-\overline{Y}_{k}\|_F^2+\beta\rho^2\|\mathcal{H}(\overline{X}_{k})-\overline{R}_{k}\|_F^2.
\end{equation}
In the other hand, we see that (\ref{lagup}) is equivalent to
\begin{eqnarray}
\sqrt{\rho}\overline{Z}_{k+1}&=&\sqrt{\rho}\overline{Z}_k+\beta\sqrt{\rho}\left(D\overline{X}_{k}-\overline{Y}_{k}\right),\\\nonumber
\sqrt{\beta}\overline{W}_{k+1}&=&\sqrt{\beta}\overline{W}_k+\rho\sqrt{\beta}\left(\mathcal{H}(\overline{X}_{k})-\overline{R}_{k}\right).
\end{eqnarray}
Using these two equalities gives 
\begin{eqnarray}
&&\left(\beta\|\overline{Z}_k\|_F^2+\rho\|\overline{W}_k\|_F^2\right)-\left(\beta\|\overline{Z}_{k+1}\|_F^2+\rho\|\overline{W}_{k+1}\|_F^2\right)\\\nonumber
&=&-2\beta\rho \left<\overline{Z}_k,D\overline{X}_{k}-\overline{Y}_{k}\right>-2\beta\rho \left<\overline{W}_k,\mathcal{H}(\overline{X}_{k})-\overline{Z}_{k}\right>-\beta^2\rho\|D\overline{X}_{k}-\overline{Y}_{k}\|_F^2-\beta\rho^2\|\mathcal{H}(\overline{X}_{k})-\overline{R}_{k}\|_F^2.
\end{eqnarray}
Using (\ref{keyf}) shows
\begin{eqnarray}
\left(\beta\|\overline{Z}_k\|_F^2+\rho\|\overline{W}_k\|_F^2\right)&-&\left(\beta\|\overline{Z}_{k+1}\|_F^2+\rho\|\overline{W}_{k+1}\|_F^2\right)\\\nonumber&\geq&-\beta\rho \left<\overline{Z}_k,D\overline{X}_{k}-\overline{Y}_{k}\right>-\beta\rho \left<\overline{W}_k,\mathcal{H}(\overline{X}_{k})-\overline{Z}_{k}\right>
\end{eqnarray}
 It follows from (\ref{keyf}) that 
\begin{equation}
\sum_{k=0}^{\infty}\left(\beta^2\rho\|D\overline{X}_{k}-\overline{Y}_{k}\|_F^2+\beta\rho^2\|\mathcal{H}(\overline{X}_{k})-\overline{R}_{k}\|_F^2\right)\leq \left(\beta\|\overline{Z}_0\|_F^2+\rho\|\overline{W}_0\|_F^2\right),
\end{equation}
which implies that $D\overline{X}_{k}-\overline{Y}_{k}\longrightarrow 0$ and $\mathcal{H}(\overline{X}_{k})-\overline{R}_{k}\longrightarrow0$ as $k\longrightarrow \infty.$\\
To show $\lim\limits_{k \rightarrow +\infty} F(R_{k}) + G(Y_{k})=F(R_*)+G(Y_*)$, we first see that the second inequality of (\ref{sadlpt}) implies 
\begin{eqnarray}
F(R_*)+G(Y_*)-F(R_k)-G(Y_k)&\leq& \left<W_*,\mathcal{H}(X_k)-R_k\right>_F+\left<Z_*,DX_k-Y_k\right>_F\\
&+&\beta\|DX_k-Y_k\|_F^2+\|\mathcal{H}(X_k)-R_k\|_F^2
\end{eqnarray}
in the other hand, by addition of (\ref{1ineq}), (\ref{2ineq}) and (\ref{3ineq}) we obtain
\begin{eqnarray}
F(R_k)+G(Y_k)-F(R_*)-G(Y_*)&\leq& -\left<W_k,\mathcal{H}(X_k)-R_k\right>_F-\left<Z_k,DX_k-Y_k\right>_F\\
&-&\beta\|DX_k-Y_k\|_F^2-\|\mathcal{H}(X_k)-R_k\|_F^2,
\end{eqnarray}
thus we have $\lim\limits_{k \rightarrow +\infty} F(R_{k}) + G(Y_{k})=F(R_*)+G(Y_*)$, i.e., objective convergence.

\section{Generalized matrix Krylov subspace for TV/L1 and TV/L2 regularizations}\label{sec4}
In this section  we will see how to generalize the generalized Krylov subspace (GKS) method proposed in \cite{LRV} to solve the generalized Sylvester matrix equation (\ref{TVL2sylv}). In \cite{LRV} GKS was introduced to solve Tikhonov regularization problems with a generalized regularization matrix. The method was next generalized in \cite{LMRS} to iteratively solve a sequence of weighted $\ell_2-$norms. It is our aim to use the fashion of the GKS method to iteratively solve the sequence of generalized Sylvester matrix equation (\ref{TVL2sylv}). Let us first introduce the following linear matrix operator
\begin{eqnarray*}
	\mathcal{A}:\mathbb{R}^{m\times n}& \rightarrow & \mathbb{R}^{m\times n} \\
	\mathcal{A}(X)&:=&A_1XA_2+A_3XA_4.
\end{eqnarray*}
the problem (\ref{TVL2sylv}) can be then expressed as follows
\begin{equation}
\mathcal{A}\left(X\right)=E_k,\quad k=0,1,...
\end{equation}
We start with the solution $X_1$ of the following linear matrix equation
\begin{equation}
\mathcal{A}\left(X\right)=E_0
\end{equation}
We search for an approximation of the solution by solving the following minimization problem,
\begin{equation}
\underset{X}{\text{ min }}\left\|\mathcal{A}\left(X\right)-E_0\right\|_F
\end{equation}
Let $X_0$ be an initial guess of $X_1$ and $P_0=\mathcal{A}\left(X\right)-E_0$ the corresponding residual. We use the modified global Arnoldi algorithm  \cite{Jb} to  construct an F-orthonormal basis $\mathcal{V}_m=\left[V_1,V_2,...,V_m\right]$ of the following matrix Krylov subspace
\begin{equation}
\mathcal{K}_m\left(\mathcal{A},P_0\right)=\text{span}\left\{P_0,\mathcal{A}\left(P_0\right),...,\mathcal{A}^{m-1}\left(P_0\right)\right\}.
\end{equation} 
This gives the following  relation
\begin{equation}
\mathcal{A}\left(\mathcal{V}_m\right)=\mathcal{V}_{m+1}\left(H_m\otimes I_n\right),
\end{equation}
where $H_m\in\mathbb{R}^{(m+1)\times m}$ is an upper Hessenberg matrix. We search for an approximated solution $X_1^m$ of $X_1$ belonging to $X_0+\mathcal{K}_m\left(\mathcal{A},P_0\right)$. This shows that $X_1^m$ can be obtained as follows
\begin{equation}\label{arnoldi}
X_1^m=X_0+\mathcal{V}_m(y_m\otimes I_n),
\end{equation}
where $y_m$ is the solution of the following reduced minimization problem
\begin{equation}
\underset{y\in\mathbb{R}^m}{\text{ min }}\left\|H_my-\|
P_0\|_Fe_1\right\|,
\end{equation}
where $e_1$ denotes the first unit vector of $R^{m+1}$.\\
Now we turn to the solutions of \begin{equation}
\mathcal{A}\left(X\right)=E_k,\quad k=1,2,...
\end{equation}
For example, in the
beginning of solving $\mathcal{A}\left(X\right)=E_1$, we reuse the F-orthonormal vectors $\mathcal{V}_m$ and we expand it to $\mathcal{V}_{m+1}=\left[\mathcal{V}_m,V_{\text{new}}\right]$, where $V_{\text{new}}$ is obtained normalizing the residual as follows
\begin{equation}
V_{\text{new}}=\frac{P_1}{\left\|P_1\right\|_F},\quad P_1=\mathcal{A}\left(X_1\right)-E_1
\end{equation}
We can then continue with $\mathcal{A}\left(X\right)=E_k$, $ k=2,3,...$ in a similar manner. Thus, at each iteration we generate the following new vector that has to be added to the generalized matrix Krylov subspace already generated to solve all the previous matrix equation,
\begin{equation}
V_{\text{new}}=\frac{P_k}{\left\|P_k\right\|_F},\quad P_k=\mathcal{A}\left(X_k\right)-E_k
\end{equation}
The idea of reusing these vectors to solve the next matrix equation, generates matrix subspaces refereed to as generalized matrix Krylov subspaces of increasing dimension \cite{Bouh}. Note that at each iteration, the residual $P_k$ is orthogonal to $\mathcal{V}_k$, since it is parallel to the gradient of the function (\ref{XTVL2}) evaluated at $X_k$.  Let $\mathcal{V}_k$ be the F-orthonormal basis of the generalized matrix Krylov subspaces at iteration $k$. When solving $\mathcal{A}\left(X\right)=E_k$, given $X_k$ and the corresponding residual $P_k$, in order to minimize the residual in the generalized matrix Krylov subspaces spanned by $\mathcal{V}_k$, we need to solve the following minimization problem
\begin{equation}\label{sylvls}
\underset{X\in\text{span}\left(\mathcal{V}_k\right)}{\text{ min }}\left\|P_k-\mathcal{A}\left(X\right)\right\|_F,
\end{equation}
The approximate solution of (\ref{sylvls}) is then given by $X_{k+1}=\mathcal{V}_k\left(y\otimes I_n\right)$. By means of the Kronecker product, we can recast (\ref{sylvls})  to a vector least-squares problem. Hence,  replacing the expression of $X_{k+1}$ into (\ref{sylvls}) yields the following minimization problem
\begin{equation}\label{vecls}
\underset{y_k}{\text{ min }}\left\|P_k-\left[\mathcal{A}(V_1),...,\mathcal{A}(V_k),\mathcal{A}(V_{\text{new}})\right]\left(y\otimes I_n\right)  \right\|_F,
\end{equation}
The problem (\ref{vecls}) can be solved by the updated version of the global QR decomposition \cite{Bou}. To use the global QR decomposition, we first need to define the $\diamond$ product. Let $A = [A_1, A_2, . . . , A_p]$ and $B = [B_1, B_2, . . . , B_\ell]$ be matrices of dimension $n\times ps$ and $n \times \ell s$,
respectively, where $A_i$ and $B_j$ $(i = 1, . . . , p; j = 1, . . . , \ell)$ are $n\times s$  matrices. Then the $p \times \ell$ matrix $A^T\diamond B$ is defined by
\begin{equation}
A^T\diamond B=
\begin{bmatrix} \left<A_1,B_1\right>_F & \left<A_1,B_2\right>_F & \cdots & \left<A_1,B_\ell\right>_F \\ 
\left<A_2,B_1\right>_F & \left<A_2,B_2\right>_F & \cdots & \left<A_1,B_\ell\right>_F\\
\vdots & \vdots & & \vdots \\
\left<A_p,B_1\right>_F & \left<A_p,B_2\right>_F & \cdots & \left<A_p,B_\ell\right>_F 
\end{bmatrix}.
\end{equation} 
Let $Q_\mathcal{A}\left(R_\mathcal{A}\otimes I_n\right)$ be the global QR of 
$\left[\mathcal{A}(V_1), \mathcal{A}(V_2), ...,\mathcal{A}(V_k)\right]$, where $Q_\mathcal{A} = [Q_1,..., Q_k]$ is an $m\times kn$ F-orthonormal matrix satisfying $Q_\mathcal{A}^T\diamond Q_\mathcal{A}=I_k$ and $R_\mathcal{A}$ is an upper triangular $k\times k$ matrix. The global QR decomposition of
$\left[\mathcal{A}(V_1),...,\mathcal{A}(V_k),\mathcal{A}(V_{\text{new}})\right]$
is defined as follows
\begin{equation}
\left[\mathcal{A}(V_1),...,\mathcal{A}(V_k),\mathcal{A}(V_{\text{new}})\right]=\left[Q_\mathcal{A},Q_{\text{new}}\right]\left(
\left[ {\begin{array}{cc}
	R_\mathcal{A} & r_\mathcal{A} \\
	0 & r_a \\
	\end{array} } \right]\otimes I_n\right),
\end{equation}
where $Q_{\text{new}}$, $r_\mathcal{A}$ and $ r_a$ are updated as follows
\begin{eqnarray}
r_\mathcal{A}&=&Q_\mathcal{A}^T\diamond\mathcal{A}(V_{\text{new}}), \quad Q=\mathcal{A}(V_{\text{new}})-Q_\mathcal{A}\left(r_\mathcal{A}\otimes I_n\right)\\ \nonumber
r_a&=&\|Q\|_F,\qquad\qquad\qquad Q_{\text{new}}=Q/r_a.
\end{eqnarray}
\begin{algorithm}
	\caption{TV/L2 for (\ref{TVL2})}\label{algo2}
	\textbf{Inputs :} $H_1$, $H_2$, $C$, $B$, $\varepsilon$\\
	\textbf{Initialization :} $X_0=B$, $Y_0=DX_0$, $Z_0=0$\\
	\textbf{Parameters :} $\mu$, $\beta$
	\begin{enumerate}
		\item Generate matrix Krylov subspace $\mathcal{V}_m$ using modified global
		Arnoldi’s process. Set $X_1=X_1^m$, where $X_1^m$ is obtained by (\ref{arnoldi})
		\item \textbf{For} $k =1,...$ until convergence, \textbf{do}
		\item Update $Y_{k}$ by (\ref{Ysol}) and $Z_{k}$ by (\ref{ZTVL2})
		\item Calculate $P_{k}=\mathcal{A}\left(X_{k}\right)-E_{k}$, where
		$E_{k}=\mathcal{H}^T(B) +D^T\left(\beta Y_{k}-Z_{k}\right)$
		\item Calculate $V_{\text{new}}=\frac{P_{k}}{\left\|P_{k}\right\|_F}$ and save $\mathcal{V}_{k+1}=\left[\mathcal{V}_k,V_{\text{new}}\right]$
		\item Update $X_{k+1}$ by solving
		$ \underset{X\in\text{span}\left(\mathcal{V}_{k+1}\right)}{\text{ min }}\left\|P_k-\mathcal{A}\left(X\right)\right\|_F$
		with the updated global QR decomposition
		\item End the iteration if 
		$\left\|X_{k+1}-X_k\right\|_F/\left\|X_k\right\|_F<\varepsilon$
	\end{enumerate}
	
\end{algorithm}
 
	  \begin{algorithm}
	  	\caption{TV/L1 for (\ref{TVL1})}\label{algo1}
 \textbf{Inputs :} $H_1$, $H_2$, $C$, $B$, $\varepsilon$\\
 \textbf{Initialization :} $R_0=\mathcal{A}\left(X_0\right)-B$, $Y_0=DX_0$, $Z_0=0$, $W_0=0$\\
	\textbf{Parameters :} $\mu$, $\beta$, $\rho$\\
		  	\begin{enumerate}
		  		\item Generate matrix Krylov subspace $\mathcal{V}_m$ using modified global
		  		Arnoldi’s process. Set $X_1=X_1^m$, where $X_1^m$ is obtained by (\ref{arnoldi})
	 \item \textbf{For} $k = ,1,...$ until convergence, \textbf{do}
	 \item Update $R_{k}$ by (\ref{RTVL1})
	and Update $Y_{k}$ by (\ref{Ysol})
	 \item Update $Z_{k}$ and $W_{k}$ by (\ref{lgup})
	 \item Calculate $P_{k}=\mathcal{A}\left(X_{k}\right)-E_{k}$, where
	 $E_{k}= \mathcal{H}^T\left(\rho R_{k}-W_{k}\right)+D^T\left(\beta Y_{k}-Z_{k}\right)$
	 \item Calculate $V_{\text{new}}=\frac{P_{k}}{\left\|P_{k}\right\|_F}$ and save $\mathcal{V}_{k+1}=\left[\mathcal{V}_k,V_{\text{new}}\right]$
	 \item Update $X_{k+1}$ by solving
	$ \underset{X\in\text{span}\left(\mathcal{V}_{k+1}\right)}{\text{ min }}\left\|P_k-\mathcal{A}\left(X\right)\right\|_F$
	with the updated global QR decomposition
	\item End the iteration if 
	$\left\|X_{k+1}-X_k\right\|_F/\left\|X_k\right\|_F<\varepsilon$
	 \end{enumerate}

	  \end{algorithm}
	  
	    \section{Numerical results}
	    This section provides some numerical results to show the performance of Algorithms TV/L1 and TV/L2
	    when applied to the restoration of blurred and noisy images. The first example applies  TV/L1 to
	    the restoration of blurred image  contaminated Gaussian blur salt-and-pepper noise while the second example apply the TV/L1 model when also a color image is contaminated by Gaussian blur salt-and-pepper noise.  The third
	    example discusses TV/L2 when applied  to the restoration of an image that have been contaminated by 
	    Gaussian blur and by additive zero-mean white
	    Gaussian noise. All computations were carried out using the MATLAB environment on an Pentium(R) Dual-Core CPU T4200 computer with 3 GB of
	    RAM. The computations were done with approximately 15 decimal digits of relative
	    accuracy. To determine the effectiveness of our solution methods, we evaluate the Signal-to-Noise
	    Ratio (SNR) defined by
	    \[\text{SNR}(X_k)=10\text{log}_{10}\frac{\|\widehat{X}-E(\widehat{X})\|_F^2}{\|X_k-\widehat{X}\|_F^2}\]
	    where $E(\widehat{X})$ denotes the mean gray-level of the uncontaminated image $\widehat{X}$. 
	    The parameters are chosen empirically to yield the best reconstruction. In all the examples we generate the matrix Krylov subspace $\mathcal{V}_1$ using only one step of the modified global Arnoldi's process.
	    \subsection*{Example 1}
\begin{figure}[ht]
	    	\begin{minipage}[b]{.45\textwidth}
	    		\centering
	    		\includegraphics[width=1\linewidth]{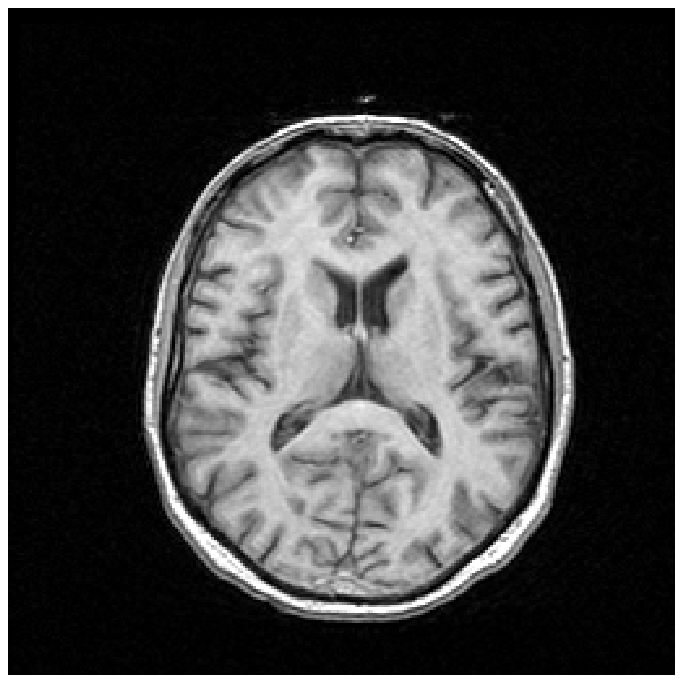}
	    		\captionof{figure}{Original image}
	    		\label{figo}
	    		\hspace{3ex}
	    	\end{minipage}%
	    	\begin{minipage}[b]{.45\linewidth}
	    		\centering
	    		\includegraphics[width=1\linewidth]{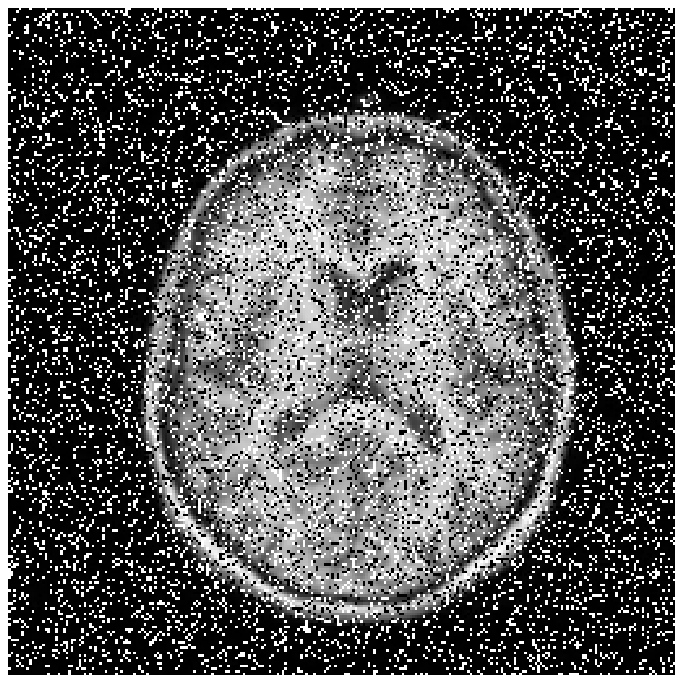}
	    		\captionof{figure}{Corrupted}
	    		\label{figbn}
	    		\hspace{1ex}
	    	\end{minipage}
	    	\begin{minipage}[b]{.45\linewidth}
	    		\centering
	    		\includegraphics[width=1\linewidth]{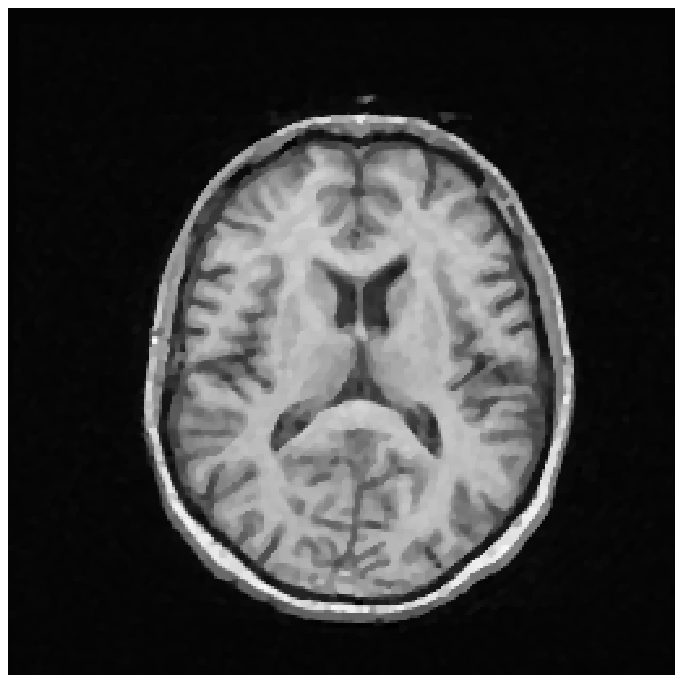}
	    		\captionof{figure}{$TV_{1}$ (SNR=19.21)}
	    		\label{figTVi}
	    		\hspace{1ex}
	    	\end{minipage}
	    	\begin{minipage}[b]{.45\linewidth}
	    		\centering
	    		\includegraphics[width=1\linewidth]{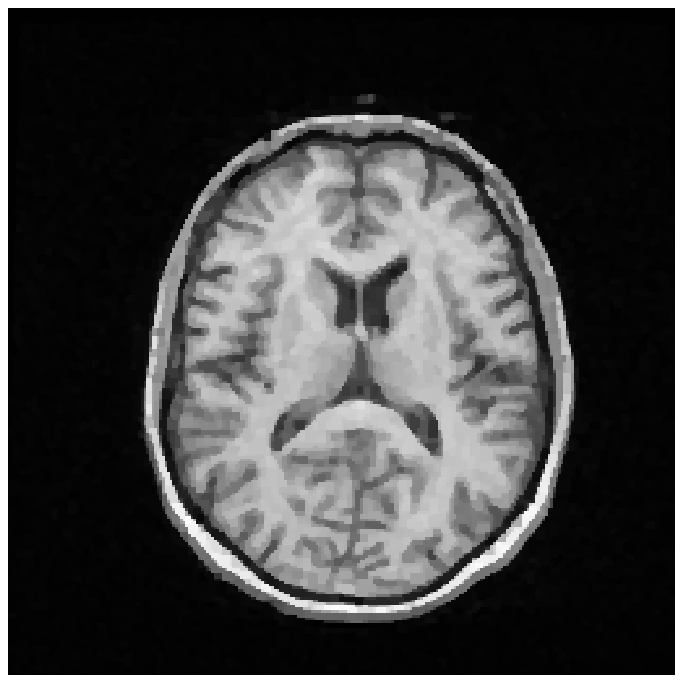}
	    		\captionof{figure}{$TV_{2}$ (SNR=17.66)}
	    		\label{figTVa}
	    		\hspace{1ex}
	    	\end{minipage}
	    \end{figure}

	    In this example the original image is the gray-scale {\tt mrin6.png} image of dimension $256\times 256$ from Matlab and it is shown in  Figure \ref{figo}. The blurring matrix $H$ is given by $H = H_1\otimes H_2\in\mathbb{R}^{256^2\times 256^2},$ where $H_1 = H_2 = [h_{ij}]$ and $[h_{ij}]$ is the Toeplitz matrix of dimension $256\times 256$ given by
	    $$ h_{ij}=
	    \left\{
	    \begin{array}{rcr}
	    \frac{1}{\sigma\sqrt{2\pi}}\text{exp}\left(-\frac{(i-j)^2}{2\sigma^2}\right), |i-j|\leq r,\\
	    0 \qquad\qquad\qquad\qquad \qquad \text{otherwise}\\
	    \end{array}
	    \right.$$
	    The blurring matrix $H$ models a blur arising in connection with the degradation of digital images by atmospheric
	    turbulence blur. We let $\sigma=1$ and $r =4$. The blurred and noisy image of Figure \ref{figbn} has been built by the product $H_2\hat{X}H_1^T$   and by adding salt-and-pepper noise of different intensity. The recovery of the image via $TV_{1}/L1$ and $TV_{2}/L1$ models  is terminated as soon as  $\left\|X_{k+1}-X_k\right\|_F/\left\|X_k\right\|_F<10^{-3}.$ Table \ref{tabtv1} report results of the performances of the $TV$ models for
	    different percentages of pixels corrupted by salt-and-pepper noise. In Figures \ref{figTVi}-\ref{figTVa}  we show the resorted images obtained
	    applying TV/L1 algorithm for $30\%$  noise level.
	    \begin{table}[h!]
	    	\centering
	    	\begin{tabular}{||c ||c|c|c||c|c| c||c|c|c|c||} 
	    		\hline
	    		&\multicolumn{3}{|c||}{Parameters}&\multicolumn{3}{|c|}{$\text{TV}_{1}$}&\multicolumn{3}{|c|}{$\text{TV}_{2}$}\\
	    		\hline
	    		Noise \% & $\mu$&$\beta$ &$\rho$& Iter &SNR&time& Iter &SNR&time  \\ [0.5ex] 
	    		\hline\hline
	    		10 & 0.05 & 50 & 5&56&23.55&10.23&141&22.64&42.55 \\ 
	    		20 & 0.1 & 50 & 5&51&21.38&8.69&106&20.16&27.16 \\
	    		30 & 0.2 & 50 & 5 &48&19.21&7.68&87&17.66&19.73\\[1ex] 
	    		\hline
	    	\end{tabular}
	    	\caption{Comparison of $TV_{1}/L1$ and $TV_{2}/L1$} models for the restoration of {\tt mrin6.png} test image corrupted by  Gaussian blur and different salt-and-pepper noise.
	    	\label{tabtv1}
	    \end{table}
	    
	    \subsection{Example 2}
	    This example illustrates the performance of TV/L1 algorithm when applied to the 
	    restoration of 3-channel RGB color images that have been contaminated by blur and salt and peppers noise.
	    The corrupted image is stored in a block vector $B$ with
	    three columns. The desired (and assumed unavailable) image is stored in the block
	    vector $\widehat{X}$ with three columns. The blur-contaminated, and noisy image
	    associated with $\widehat{X}$, is stored in the block vector $B$. 
	    
\begin{figure}
	    	\begin{minipage}{.45\linewidth}
	    		\centering
	    		\includegraphics[width=1\linewidth]{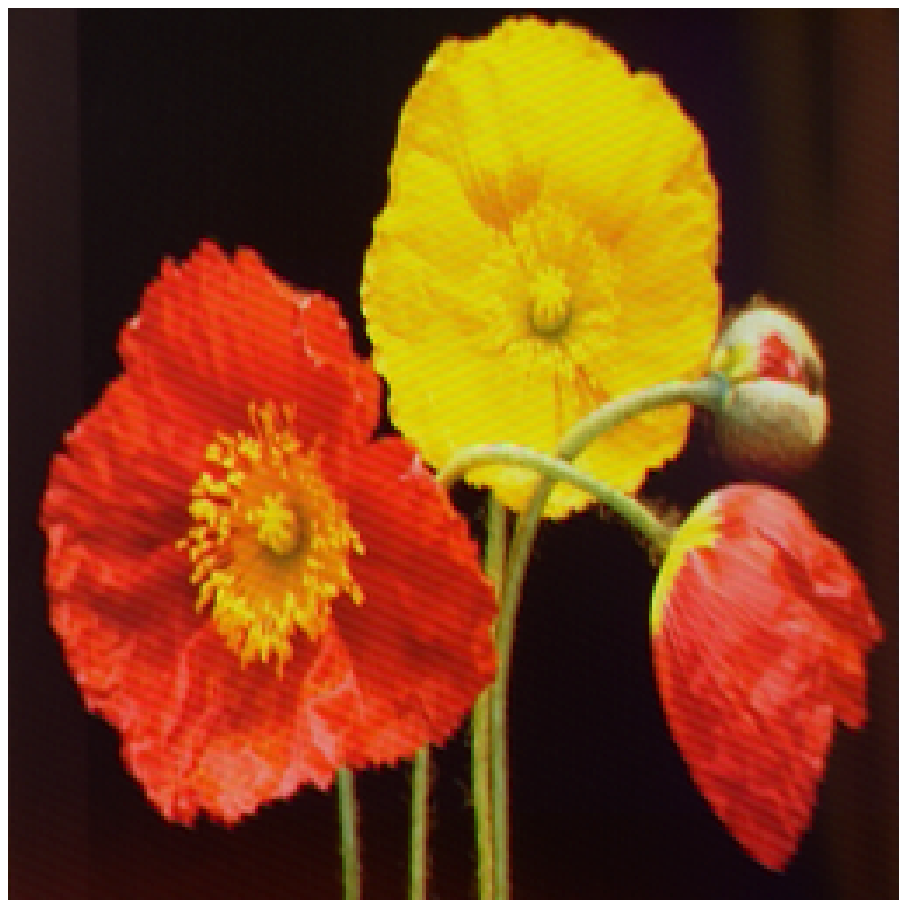}
	    		\captionof{figure}{Original image}
	    		\label{papavo}
	    	\end{minipage}%
	    	\begin{minipage}{.45\textwidth}
	    		\centering
	    		\includegraphics[width=1\linewidth]{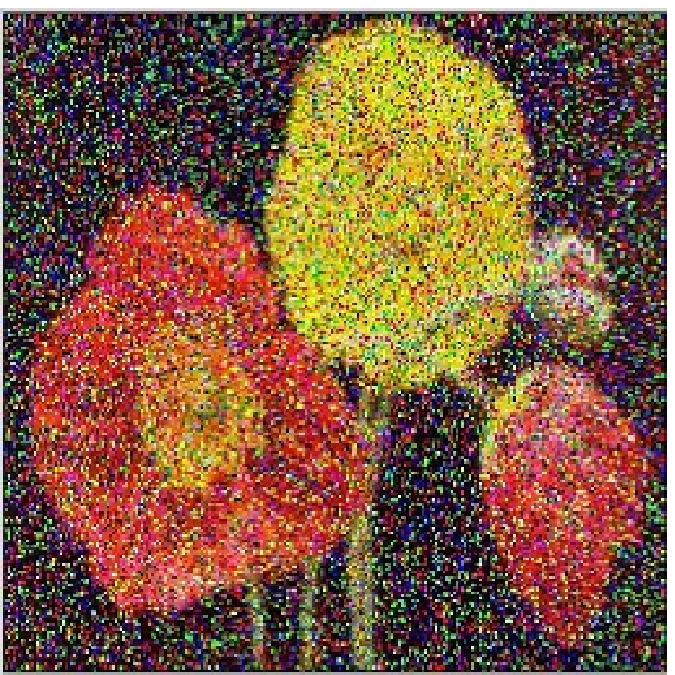}
	    		\captionof{figure}{Corrupted}
	    		\label{papavbn}
	    	\end{minipage}
	    	\begin{minipage}{.45\linewidth}
	    		\centering
	    		\includegraphics[width=1\linewidth]{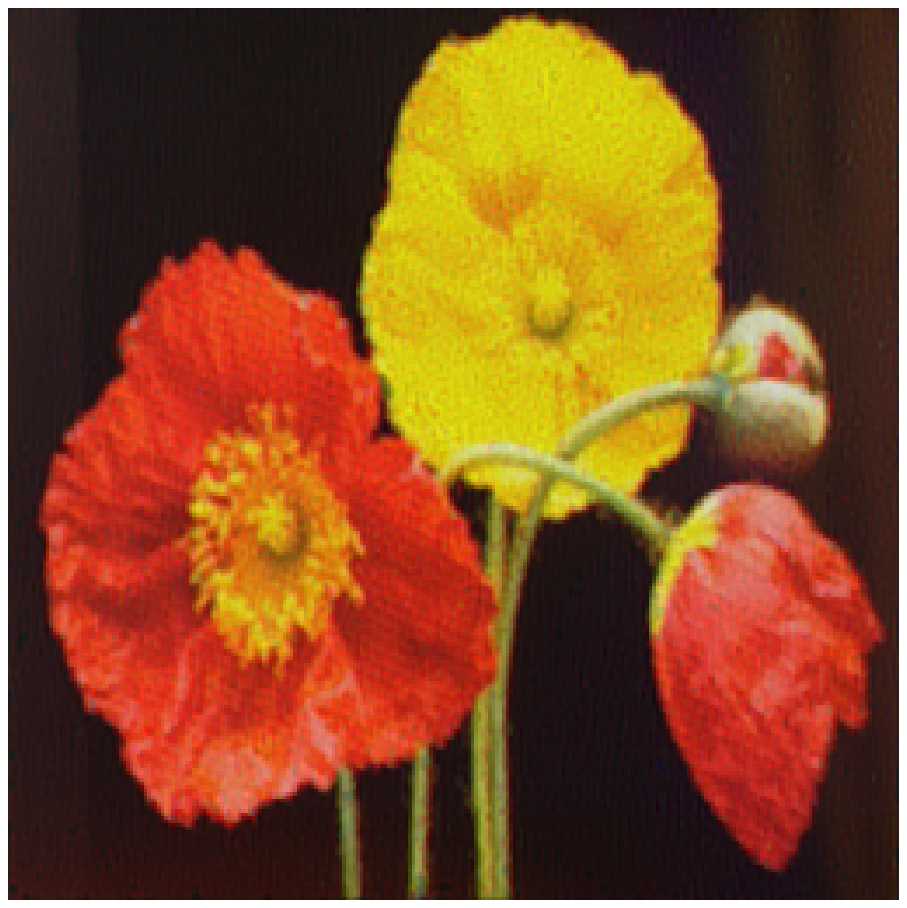}
	    		\captionof{figure}{$TV_{1}$ (SNR=20.90)}
	    		\label{papavTVi}
	    	\end{minipage}
	    	\begin{minipage}{.45\linewidth}
	    		\centering
	    		\includegraphics[width=1\linewidth]{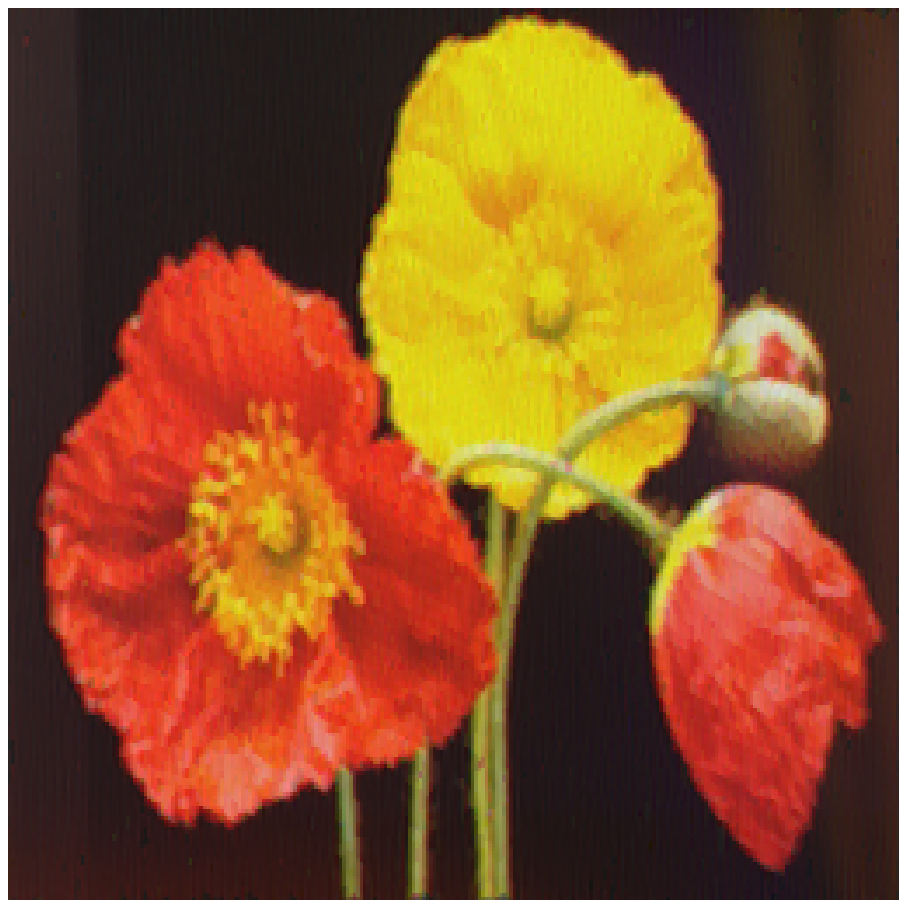}
	    		\captionof{figure}{ $TV_{2}$ (SNR=21.13)}
	    		\label{papavTVa}
	    	\end{minipage}
	    \end{figure}

	    We consider the within-channel blurring only. Hence the blurring matrix $H_1$ in
	    \eqref{model1} is the $3\times3$ identity matrix. The blurring matrix $H_2$ in \eqref{model1},
	    which describes the blurring within each channel, models Gaussian blur and is determined  with the MATLAB function {\sf blur} from \cite{Regtools}. This
	    function has two parameters, the half-bandwidth of the Toeplitz blocks $r$ and the 
	    variance $\sigma$ of the Gaussian PSF. For this example we let $\sigma=1$ and $r=4$. The
	    original (unknown) RGB image $\widehat{X}\in 256\times 256\times 3$ is the $\tt papav256$
	    image from MATLAB. It is shown in Figure \ref{papavo}. The 
	    associated blurred and noisy image $B$ with $30\%$ noise level is shown in Figure \ref{papavbn}. Given the contaminated image $B$, we 
	    would like to recover an approximation of the original image $\widehat{X}$. The recovery of the image via $TV_{1}/L1$ and $TV_{2}/L1$ models  is terminated as soon as  $\left\|X_{k+1}-X_k\right\|_F/\left\|X_k\right\|_F<10^{-2}.$ Table
	    \ref{tablergb} compares the results obtained by $TV_{1}/L1$ and $TV_{2}/L1$ models. 
	    
	    The restorations obtained with
	    $TV_{1}/L1$  and $TV_{2}/L1$ for noise level $30\%$ are shown in Figure 
	    \ref{papavTVi} and the Figure 
	    \ref{papavTVa}, respectively. 
	    \begin{table}[h!]
	    	\centering
	    	\begin{tabular}{||c ||c|c|c||c|c| c||c|c|c|c||} 
	    		\hline
	    		&\multicolumn{3}{|c||}{Parameters}&\multicolumn{3}{|c|}{$\text{TV}_{1}$}&\multicolumn{3}{|c|}{$\text{TV}_{2}$}\\
	    		\hline
	    		Noise \% & $\mu$&$\beta$ &$\rho$& Iter &SNR&time& Iter &SNR&time  \\ [0.5ex] 
	    		\hline\hline
	    		10 & 0.1 & 80 & 5&13&24.66&9.01&14&24.32&9.73 \\ 
	    		20 & 0.125 & 80 & 5&17&23.00&12.64&17&22.71&12.36 \\
	    		30 & 0.125 & 80 & 5 &19&20.90&13.35&19&21.13&13.89\\[1ex] 
	    		\hline
	    		
	    	\end{tabular}
	    	\caption{Comparison of $TV_{1}/L1$ and $TV_{2}/L1$} models for the restoration of {\tt papav256.png} test colour image corrupted by  Gaussian blur and different salt-and-pepper noise.
	    	\label{tablergb}
	    \end{table}
	    \subsection{Example 3}
\begin{figure}[ht]
	    	\begin{minipage}{0.45\linewidth}
	    		\centering
	    		\includegraphics[width=1\linewidth]{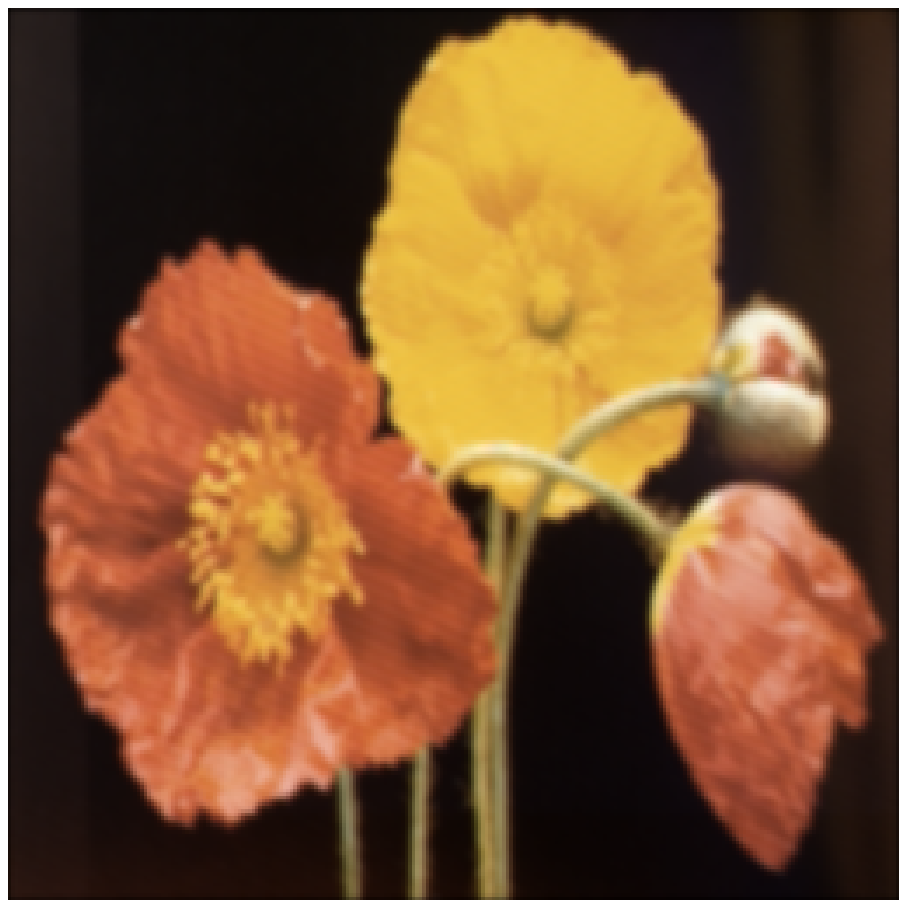}
	    		\captionof{figure}{Blurred image}
	    		\label{papav256}
	    	\end{minipage}%
	    	\begin{minipage}{.45\linewidth}
	    		\centering
	    		\includegraphics[width=1\linewidth]{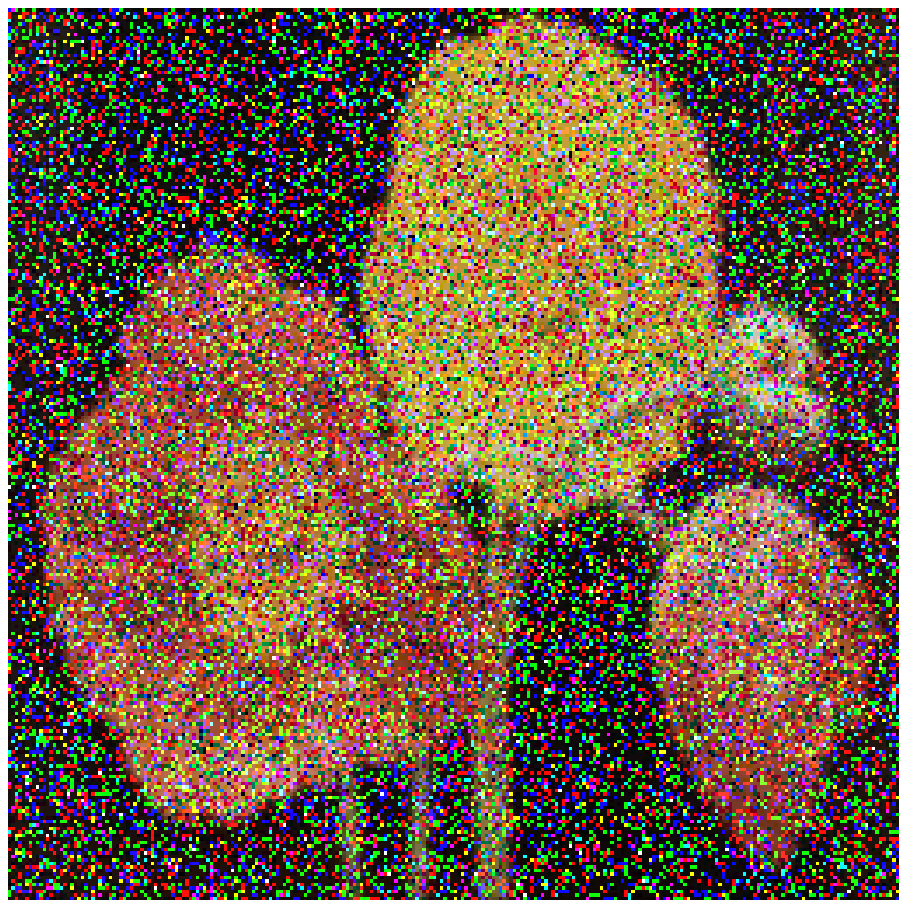}
	    		\captionof{figure}{Blurred and noisy image}
	    		\label{papav256bn}
	    	\end{minipage}
	    	\begin{minipage}{.45\linewidth}
	    		\centering
	    		\includegraphics[width=1\linewidth]{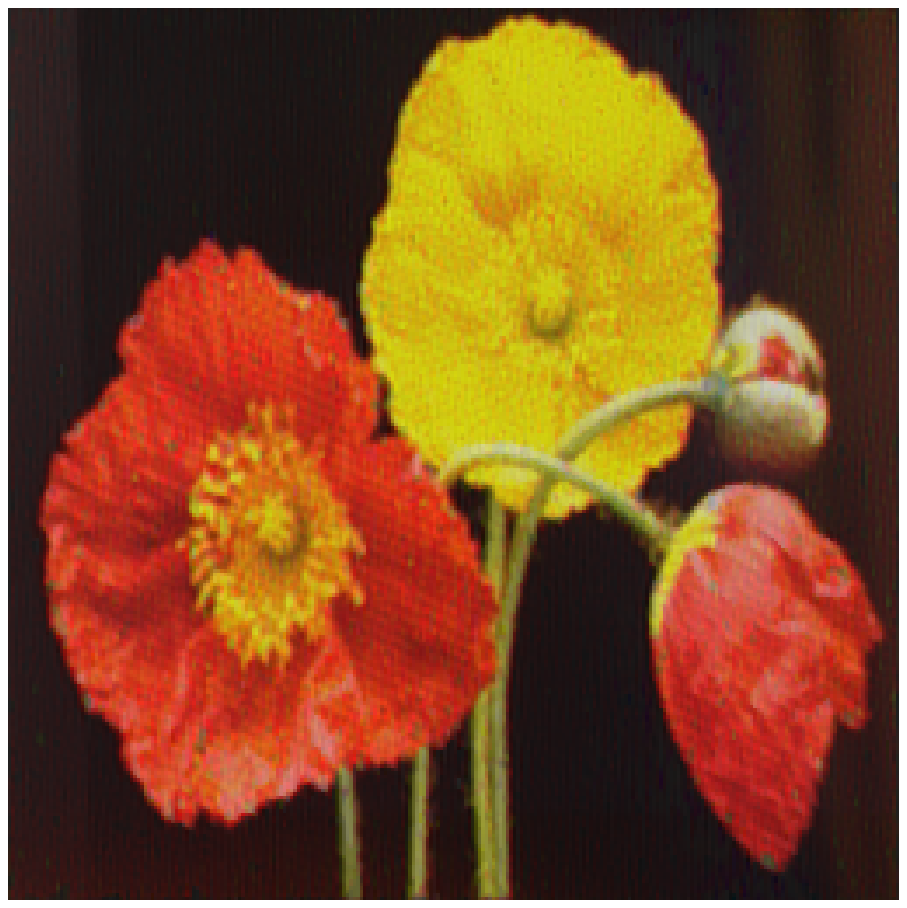}
	    		\captionof{figure}{$TV_{1}$ (SNR=19.50)}
	    		\label{papav256TVa}
	    	\end{minipage}
	    	\begin{minipage}{.45\linewidth}
	    		\centering
	    		\includegraphics[width=1\linewidth]{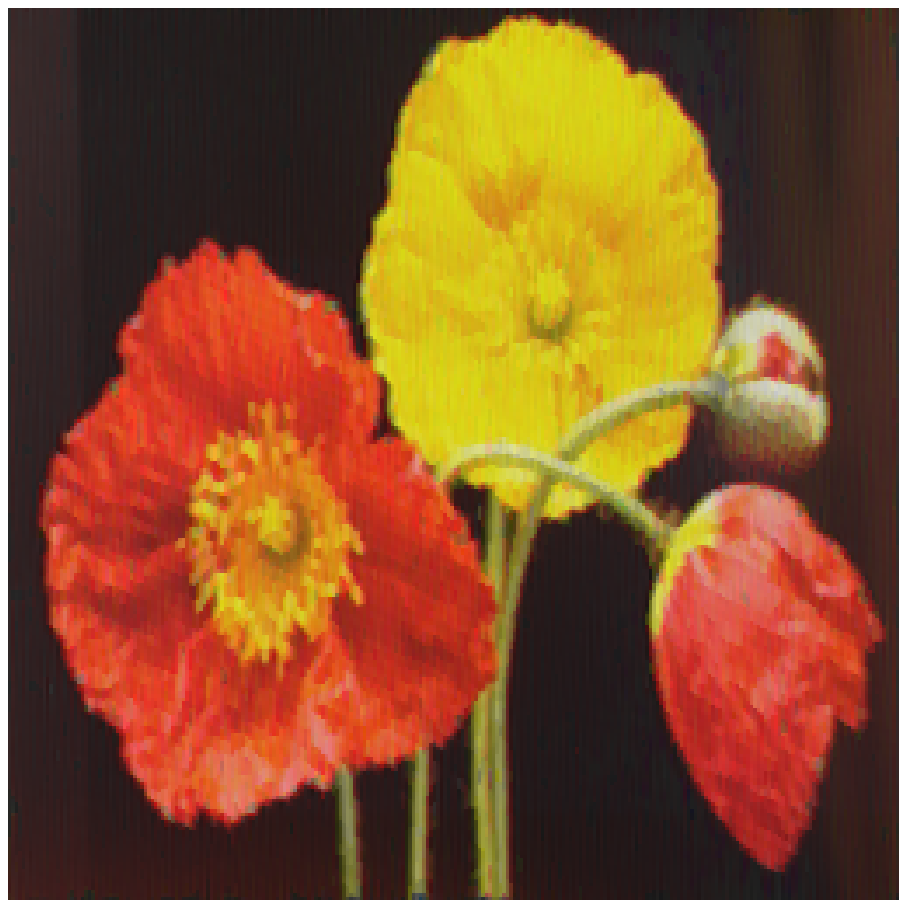}
	    		\captionof{figure}{ $TV_{2}$ (SNR=19.90)}
	    		\label{papav256TVi}
	    	\end{minipage}
	    \end{figure}

	    In this example we present the experimental results recovered by Algorithm 1  for the reconstruction of a cross-channel blurred image. We consider the same original RGB image  and  the same within-channel blurring matrix $H_1$, as in Example 2, with the same parameters. The cross-channel blurring  is determined by a matrix $H_2$. In our example we let $H_2$ to be
	    \[
	    H_2=\begin{bmatrix}
	    0.7&
	    0.2&0.1\\0.25&0.5&0.25\\0.15&0.1&0.75
	    \end{bmatrix}.
	    \]
	    This matrix is obtained from \cite{HNO}. The cross-channel blurred  image without noise is represented by $H_1\widehat{X}H_2^T$ and it is shown in Figure (\ref{papav256}) . The 
	    associated blurred and noisy image $B$ with $30\%$ noise level is shown in Figure (\ref{papav256bn}).  The cross-channel blurred  and noisy image has been reconstructed using Algorithm 1 as soon as $\left\|X_{k+1}-X_k\right\|_F/\left\|X_k\right\|_F<10^{-2}.$ The  restored images obtained with  TV/L1 models are shown in  Figures (\ref{papav256TVi})-(\ref{papav256TVa}). 
	    \subsection{Example 4  }	
\begin{figure}[ht]
	    	\begin{minipage}{.45\linewidth}
	    		\centering
	    		\includegraphics[width=1\linewidth]{mrin6.eps}
	    		\captionof{figure}{Original image}
	    		\label{mrino}
	    	\end{minipage}%
	    	\begin{minipage}{.45\linewidth}
	    		\centering
	    		\includegraphics[width=1\linewidth]{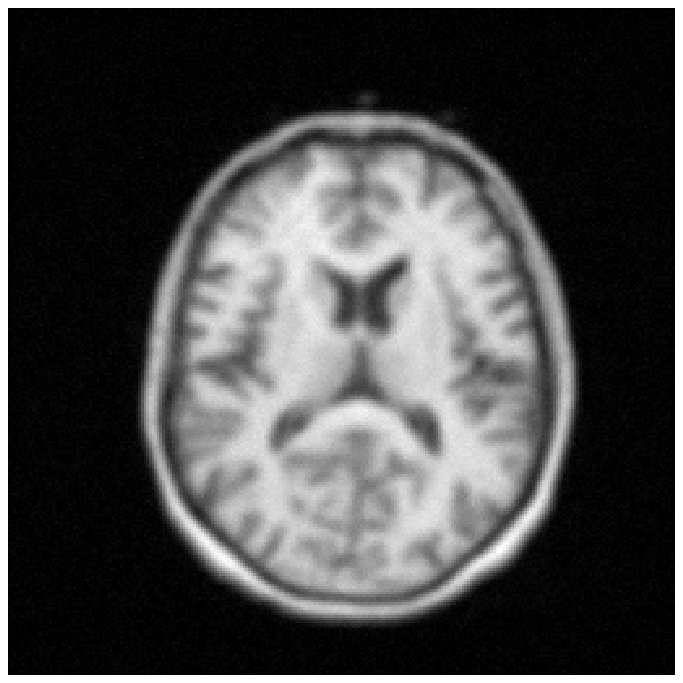}
	    		\captionof{figure}{Corrupted}
	    		\label{mrinbn}
	    	\end{minipage}
	    	\begin{minipage}{.45\linewidth}
	    		\centering
	    		\includegraphics[width=1\linewidth]{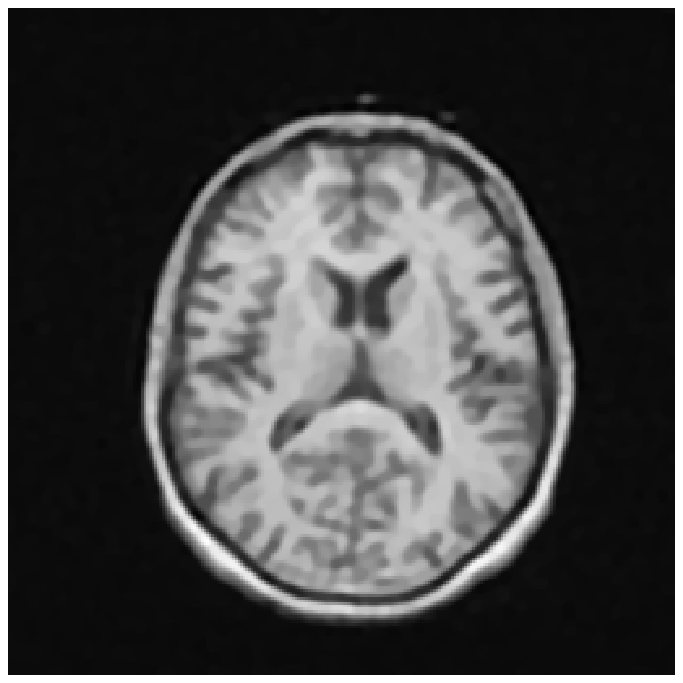}
	    		\captionof{figure}{$TV_{1}$ (SNR=15.70)}
	    		\label{mrinTVi}
	    	\end{minipage}
	    	\begin{minipage}{.45\linewidth}
	    		\centering
	    		\includegraphics[width=1\linewidth]{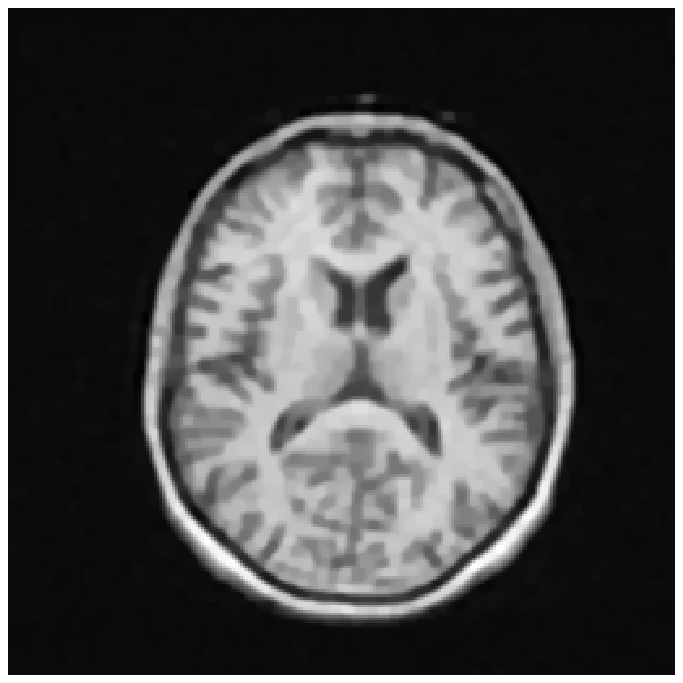}
	    		\captionof{figure}{ $TV_{2}$ (SNR=15.60)}
	    		\label{mrinTVa}
	    	\end{minipage}
	    \end{figure}

	    In this example we consider the restoration of the gray-scale {\tt mrin6.png} image degraded by the same blurring matrices $H_1$ and $H_2$ defined in Example 1 with $\sigma=2$ and $r=4$, and by additive zero-mean white
	    Gaussian noise with different different noise levels. This noise level  is  defined as follows 
	    $\nu=\frac{||E||_F}{||\widehat{B}||_F}
	    $,
	    where  $E$ denotes the block
	    vector that represents the noise in $B$, i.e., $B:=\widehat{B}+E$, and $\widehat{B}$ is the noise-free image
	    associated with original image  $\widehat{X}$. For this kind of noise, we consider the $TV_{1}/L2$ and $TV_{2}/L2$ models.  The recovery of the image via $TV_{1}/L1$ and $TV_{2}/L1$ models  is terminated as soon as  $\left\|X_{k+1}-X_k\right\|_F/\left\|X_k\right\|_F<10^{-3}.$
	    In Table \ref{tablemrin}, we compare the results obtained by $TV_{1}/L2$ and $TV_{2}/L2$ for different noise levels. Figure \ref{mrinbn} shows the image degraded by $0.01$ noise level. Figure \ref{mrinTVi} and Figure \ref{mrinTVa} show the restored images obtained by $TV_{1}/L2$ and $TV_{2}/L2$, respectively.
	    \begin{table}[h!]
	    	\centering
	    	\begin{tabular}{||c ||c|c||c|c| c||c|c|c|c||} 
	    		\hline
	    		&\multicolumn{2}{|c||}{Parameters}&\multicolumn{3}{|c|}{$\text{TV}_{1}$}&\multicolumn{3}{|c|}{$\text{TV}_{2}$}\\
	    		\hline
	    		Noise \% & $\mu$&$\beta$ & Iter &SNR&time& Iter &SNR&time  \\ [0.5ex] 
	    		\hline\hline
	    		0.001 & 0.0001 & 0.1 & 53&18.32&9.30&52&18.32&10.10 \\ 
	    		0.01 & 0.001 & 30&20&15.70&2.65&21&15.60&2.60\\[1ex] 
	    		\hline
	    	\end{tabular}
	    	\caption{Comparison of $TV_{1}/L2$ and $TV_{2}/L2$} models for the restoration of {\tt imrin6.png} test  image corrupted by  Gaussian blur and different white Gaussian noise level.
	    	\label{tablemrin}
	    \end{table} 
	    \subsection{Example 5}
\begin{figure}[ht]
	     	\begin{minipage}{0.45\linewidth}
	     		\centering
	     		\includegraphics[width=1\linewidth]{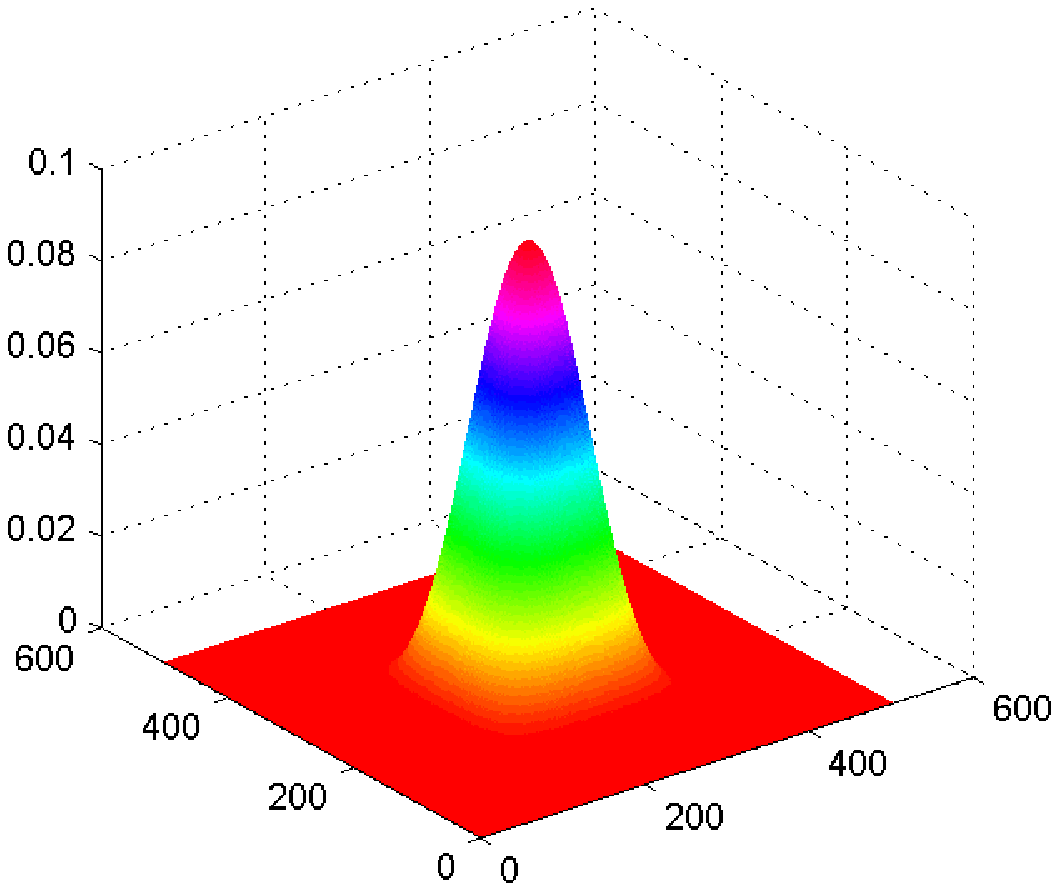}
	     		\captionof{figure}{True object}
	     		\label{phillipsor}
	     	\end{minipage}%
	     	\begin{minipage}{.45\linewidth}
	     		\centering
	     		\includegraphics[width=1\linewidth]{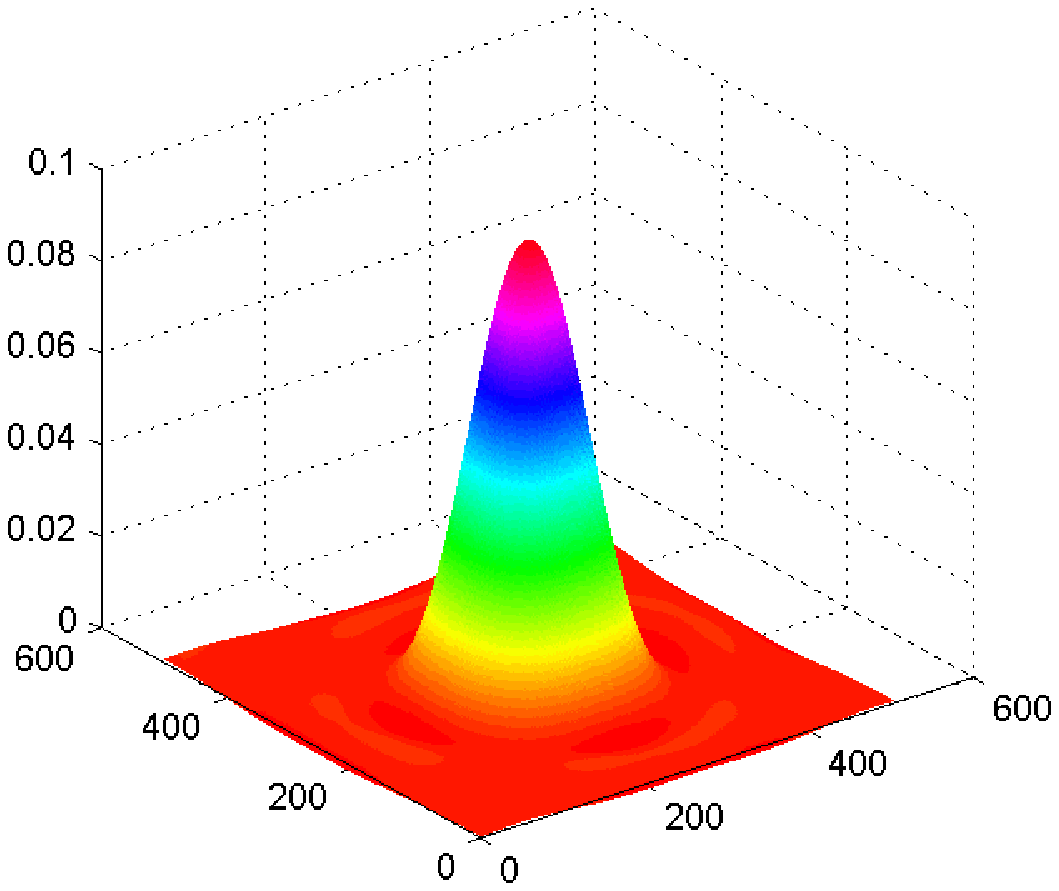}
	     		\captionof{figure}{Approximate solution}
	     		\label{phillipsrest}
	     	\end{minipage}
	     \end{figure}

	    In this example, we consider the Fredholm integral equation
	    \begin{equation}\label{fredholm}
	    \int\int_\Omega K(x,y,s,t)f(s,t)dsdt=g(x,y),\qquad (x,y)\in\Omega,
	    \end{equation}
	    where $\Omega=[-6,6]\times[-6,6]$. Its kernel, solution,  and right-hand side are given by
	    $$K(x,y,s,t)=k_1(x,s)k_1(y,t),\quad (x,y)\in\Omega,\quad (s,t)\in\Omega,$$
	    $$f(x,y)=f_1(x)f_1(y),$$ 
	    $$g(x,y)=g_1(x)g_1(y),$$
	    where
	    $$f_1(s):=\left\{
	    \begin{array}{cl}
	    1+\cos(\frac{\pi}{3}s),& \quad  |s|\leq \frac{\pi}{3},\\
	    0,& \quad \text{otherwise}.
	   \end{array}
	   \right.$$
	    $$k_1(s,x):=f_1(s-x)$$ $$g_1(s):=(6-|s|)\left(1+\frac{1}{2}\cos\left(\frac{\pi}{3}s\right)\right)+\frac{9}{2\pi}\sin\left(\frac{\pi}{3}|s|\right).$$
	    We use the code {\tt phillips} from Regularization Tools \cite{Regtools} to discretize (\ref{fredholm}) by a Galerkin
	    method with orthonormal box functions as test and trial functions to obtain $H_1$ and $H_2$ of size $500$. From the output of the code {\tt phillips} we determine a scaled approximation $\widehat{X}\in \mathbb{R}^{500\times 500}$ of the exact solution $f(x,y)$. Figure \ref{phillipsor} displays this exact solution. To determine the effectiveness of our approach, we evaluate the relative error 
	    $$\text{Re}=\frac{||\widehat{X}-X_k||_F}{||\widehat{X}||_F}$$
	    of the computed approximate solution $X_k$ 
	    obtained with Algorithm 1. Table \ref{tabphillips} shows the relative error in
	    approximate solutions determined by Algorithm 1 for different noise
	    levels, as well as the number of iterations required to satisfy $\left\|X_{k+1}-X_k\right\|_F/\left\|X_k\right\|_F<10^{-3}.$
	    Figure \ref{phillipsrest} displays the computed approximate solution obtained when the noise
	    level is $0.1$.
	     \begin{table}[h!]
	     	\centering
	     	\begin{tabular}{||c ||c|c||c|c| c||c|c|c|c||} 
	     		\hline
	     		&\multicolumn{2}{|c||}{Parameters}&\multicolumn{3}{|c|}{$\text{TV}_{1}$}&\multicolumn{3}{|c|}{$\text{TV}_{2}$}\\
	     		\hline
	     		Noise \% & $\mu$&$\beta$ & Iter &Re&time& Iter &Re&time  \\ [0.5ex] 
	     		\hline\hline
	     		0.001 & 0.0001 & 0.1 & 12&$4.01\times10^{-2}$&9.05&9&$4.71\times10^{-2}$&6.52 \\ 
	     		0.01 & 0.001 & 30&13&$3.99\times10^{-2}$&9.63&13&$3.98\times10^{-2}$&9.66\\ 
	     		0.1 & 0.1 & 40&15&$4.07\times10^{-2}$&10.94&15&$4.07\times10^{-2}$&11.38\\[1ex] 
	     		\hline
	     	\end{tabular}
	     	\caption{Comparison of $TV_{1}/L2$ and $TV_{2}/L2$} models for the solution of (\ref{fredholm}) with different white Gaussian noise level.
	     	\label{tabphillips}
	     \end{table} 
	  \bibliographystyle{model1b-num-names}
	  
	\end{document}